\documentclass[10pt]{amsart}

\usepackage{verbatim}

\usepackage{amscd,amssymb,hyperref}
\usepackage[mathscr]{eucal}
\usepackage{graphicx}

\theoremstyle{plain}

\numberwithin{equation}{section}

\newtheorem{theorem}{Theorem}[section]
\newtheorem{lemma}[theorem]{Lemma}
\newtheorem{proposition}[theorem]{Proposition}
\newtheorem{corollary}[theorem]{Corollary}

\theoremstyle{plain}
\newtheorem{claim*}{Claim}

\theoremstyle{definition}

\newtheorem{definition}{Definition}[section]

\theoremstyle{remark}
\newtheorem{remark}{Remark}[section]
\newtheorem*{example*}{Example}
\newtheorem*{remark*}{Remark}

\newcommand{\M}{\mathcal{M}}

\DeclareMathOperator{\Tr}{Tr}
\let\Re=\undefined
\DeclareMathOperator{\Re}{Re}

\DeclareMathOperator{\vol}{vol}

\DeclareMathOperator{\Diff}{\textit{Diff}}

\newcommand{\mat}{ \begin{pmatrix} }
\newcommand{\emat}{\end{pmatrix}}

\DeclareMathOperator{\dist}{dist}

\begin{document}

\title[Uniformizations, Determinants of Laplacians, Isospectrality]
{Surfaces with boundary: their uniformizations, determinants of
Laplacians, and isospectrality}

\author{Young-Heon Kim}

\address{Department of Mathematics, University of Toronto,
Toronto, Ontario, Canada M5S 2E4}

\email{yhkim@math.toronto.edu}

\date{\today}


\subjclass[2000]{58J53; 32G15}

\begin{abstract}
  Let $\Sigma$ be a compact surface of type $(g, n)$, $n > 0$,
  obtained by removing $n$ disjoint
  disks from a closed surface of
  genus $g$. Assuming $\chi(\Sigma)<0$, we show that on $\Sigma$, the set of flat metrics which have
  the same Laplacian spectrum of
  Dirichlet boundary condition
  is compact in the $C^\infty$ topology.
  This isospectral compactness extends the result of Osgood, Phillips,
  and Sarnak \cite{OPS3} for type $(0,n)$  surfaces, whose examples
  include bounded plane domains.

  Our main ingredients are as following. We first show that the determinant of the
  Laplacian is a proper function
  on the moduli space of geodesically bordered hyperbolic
  metrics on $\Sigma$. Secondly, we show that the space of such metrics
  is homeomorphic (in the $C^\infty$-topology) to the space
  of flat metrics (on $\Sigma$) with constantly curved boundary.
  Because of this, we next reduce the complicated degenerations
  of flat metrics to the simpler and well-known
  degenerations of hyperbolic metrics, and
  we show that determinants of Laplacians of flat
  metrics on $\Sigma$, with fixed area and
  boundary of constant geodesic curvature, give
  a proper function on the corresponding moduli space.
  This is interesting
  because
  Khuri \cite{Kh} showed that if the boundary length (instead of the area) is fixed, the determinant is not a proper function
  when $\Sigma$ is of type $(g, n) , \ g>0$;
  while Osgood, Phillips, and Sarnak \cite{OPS3} showed the
  properness when $g=0$.
\end{abstract}

\maketitle


\section{Introduction}
Kac's \cite{Ka} famous question, `Can one hear the shape of a
drum?' asks whether we can determine a Riemannian
manifold by knowing its Laplacian spectrum.  Although some
Riemannian manifolds are determined uniquely by their spectra
(for example, the two dimensional round sphere),
 there are many
counter-examples, in particular, continuous
families of Riemannian metrics on some compact manifolds which are
isospectral but not (locally) isometric. (See Gordon's survey article
\cite{Go}.) Thus it is important to know  the size of the set
of all the Riemannian metrics  with the same spectrum on a given compact manifold.
This paper addresses the question of whether this isospectral set is compact
in the $C^\infty$-topology. A sequence $\{ \sigma_i \}$ of
Riemannian metrics on a compact manifold $M$ is said to converge to
a Riemannian metric $\sigma$ on $M$ in the $C^\infty$-topology if there
exist diffeomorphisms $F_i$ of $M$ such that $F_i ^* \sigma_i$
converge to $\sigma$ in the $C^\infty$ sense. In particular, metrics
in a compact set in the $C^\infty$-topology are all quasi-isometric by
uniform constants.

This paper focuses on Riemannian metrics on compact orientable bordered surfaces
and their Laplacians (denoted $\Delta$) on functions with Dirichlet
boundary condition.  The following is the
first main conclusion.

\begin{theorem}\label{compact isospectral}
  (See
  Theorem~\ref{T:isospectral compactness}.)
  Let $\Sigma$ be a compact orientable surface with boundary and assume
  the Euler characteristic $\chi(\Sigma) < 0$. The set of all the smooth flat
  (i.e. zero Gaussian curvature) metrics on $\Sigma$
  which have the same Laplacian spectrum of Dirichlet boundary condition is compact in
  the $C^\infty$-topology.
 \end{theorem}

Theorem~\ref{compact isospectral} extends one of the famous results
of Osgood, Phillips, and Sarnak, who showed
$C^\infty$-compactness for isospectral sets of  bounded plane domains
\cite{OPS3}. (They also showed
$C^\infty$-compactness of isospectral sets \cite{OPS2} for closed surfaces.)
Theorem~\ref{compact isospectral} allows us to deal with new examples such as 
the flat surfaces obtained by removing arbitrary
(smooth) neighborhoods of each vertices from compact 2-dimensional
simplicial complexes which are manifolds. The topologies of these examples
can be much more complicated
than those of plane domains.

As in the works of Osgood, Phillips, and Sarnak,  the
\emph{determinant $\det \Delta$ of the Laplacian} is used as our main
analytical notion. It was first introduced by Ray and Singer
\cite{RS}\cite{RS2} and has become one of the central objects in geometric analysis,
algebraic geometry, and string theory. It is defined (see
Section~\ref{Heights}) by using the analytic continuation of the
spectral zeta function
\begin{align*}
  \zeta(s) = \sum_{0 < \lambda  \in \text{Spec}(\Delta)} \lambda^{-s}
\end{align*}
to the origin and the formula
\begin{align*}
  -\log \det \Delta = \zeta ' (0) .
\end{align*}

Osgood, Phillips, and Sarnak call $-\log \det \Delta$ the
\emph{height} of the Riemannian metric. They use the
height as a function on the moduli space of Riemannian metrics to
study isospectral problems, given the obvious fact that isospectral
metrics have the same height.

Osgood, Phillips, and Sarnak first analyzed the extremal properties of the height function $h$ in each
conformal class of a surface and showed the uniformization theorem \cite{OPS1}.
Namely,
in each conformal class of Riemannian metrics on a compact surface, if
there is no boundary, there is a unique metric of constant curvature;
if the surface has boundary, there is a unique \emph{uniform metric of type
I}, i.e. a constant
curvature metric with geodesic boundary, and
a unique \emph{uniform metric of type II}, i.e. a flat metric with
constant geodesic curvature boundary.
These uniform metrics realize the minimum of the height under
certain constraints.

The above uniformization theorem allows Osgood, Phillips, and Sarnak
to reduce the isospectral compactness problem to the properness of the height function
on the moduli space of uniform metrics \cite{OPS2} \cite{OPS3}. This properness was proved by Wolpert \cite{W}
(also by Bismut and Bost \cite{BB} in algebraic geometry context) for closed hyperbolic surfaces,
and by Osgood, Phillips, and Sarnak on the moduli spaces of uniform metrics of type II with
fixed boundary length on punctured spheres \cite{OPS3}.
The proof of Osgood et al. is quite involved, mainly due to the
complicated degeneration patterns of flat metrics, in comparison to cases
involving hyperbolic surfaces (or uniform metrics of type I for
nonempty boundary case) where degenerations occur simply when one
pinches closed geodesics (\emph{thick-thin decomposition}).
Moreover, this properness due to Osgood et al. cannot be extended to the higher genus
case as Khuri \cite{Kh} showed that the
height is \emph{not} a proper function on the moduli space of
uniform metrics of type II with fixed boundary length when the
base surface is of type $(g, n) , \ g>0$, contrasting to the case of type $(0, n)$ surfaces in \cite{OPS3}. A surface of type $(g,n)$ is the surface obtained by
removing $n$ disjoint disks from a closed surface of genus $g$.

Our idea is to use the analysis of hyperbolic side (uniform metrics
of type I) to get results for the flat side (uniform metrics of
type II). To describe the key lemma for this connection,
first denote $\mathbf{M}_I (\Sigma, A)$ and
$\mathbf{M}_{II} (\Sigma, A)$
(see Definition~\ref{space of uniform metrics}) as the space of uniform metrics of type
I and type II, respectively, with fixed area $A$ on a compact orientable
surface $\Sigma$ with boundary.  These spaces induce the corresponding moduli spaces
$\mathcal{M}_I (\Sigma, A)$, $\mathcal{M}_{II}
(\Sigma, A)$, after taking quotient by the group $\Diff \, (\Sigma)$
of diffeomorphisms of $\Sigma$. Then the following
theorem is proved.
\begin{theorem}\label{two uniformization homeo}(See Theorem~\ref{Psi})
  Let $\Sigma$ be a compact orientable surface with boundary
  and assume $\chi(\Sigma) <0$.
  The two spaces
  $\mathbf{M}_I (\Sigma, A)$ and $\mathbf{M}_{II}(\Sigma, A)$
  are homeomorphic
  in the $C^\infty$-topology, and
  so are $\mathcal{M}_I (\Sigma , A)$ and $\mathcal{M}_{II}(\Sigma, A)$.
\end{theorem}

An important fact used in the hyperbolic (or type I) case is the following theorem concerning
the properness of the height of hyperbolic surfaces with geodesic
boundary.

\begin{theorem}\label{proper on hyperbolic moduli intro}
  (See Corollary~\ref{proper on hyperbolic moduli}
  or Theorem~\ref{T:proper on moduli I})
  Let $\mathcal{M}$ be the moduli space
  of compact hyperbolic
  (Gaussian curvature $\equiv -1$)
  surfaces with geodesic boundary.
  The height function $h$ on $\mathcal{M}$ is proper, i.e.
  \begin{align*}
    h(M) \rightarrow +\infty
  \end{align*}
  as the isometry class $[M]$ approaches $\partial \mathcal{M}$.
\end{theorem}
This theorem is a corollary of an asymptotic inequality for the height (see
Theorem~\ref{hyperbolic height}) which we obtain
using the so called \emph{insertion lemma}. This lemma, first introduced by Sarnak \cite{Sa2},
uses thick-thin decompositions of hyperbolic surfaces. Our
asymptotic inequality partially
extends the asymptotic formula of Wolpert \cite{W} or of
Bismut and Bost \cite{BB} (see also \cite{Lu}).

To prove Theorem~\ref{compact isospectral} we use Theorem~\ref{two
uniformization homeo}, Theorem~\ref{proper on hyperbolic moduli
intro}, and a method of Osgood, Phillips, and Sarnak \cite{OPS3},
which uses the work by Melrose \cite{Me} for boundary geodesic
curvature of isospectral flat surfaces.
Our approach gives both an extension of the result in \cite{OPS3}
and a simpler treatment.

On the height of uniform metrics of type II (or
 flat metrics with boundary of constant geodesic curvature),
we easily get the following theorem from
Theorem~\ref{two uniformization homeo} and
Theorem~\ref{proper on hyperbolic moduli intro}.
\begin{theorem}\label{proper on moduli II intro}
  (See Theorem~\ref{T:proper on moduli II}.)
  Suppose
  $\chi(\Sigma) <0 $. For each $A >0$,
  the height $h$ is a proper function on the moduli space
  $ \mathcal{M}_{II} (\Sigma, A)$, i.e.
  \begin{align*}
    h(M) \rightarrow + \infty
  \end{align*}
  as the isometry class $[M]$ approaches
  $ \partial \mathcal{M}_{II}(\Sigma, A) $.
\end{theorem}
This result is remarkable in comparison with the above non-properness
result of Khuri \cite{Kh}. It should be interesting to see
the reason why the two conditions, one fixes the boundary
length and the other fixes the area, result so differently in the
heights.

Notice that the properness in
Theorem~\ref{proper on hyperbolic moduli intro}
and Theorem~\ref{proper on moduli II intro}
have another interesting feature: they insure
the existence of the global minimum. On the variational
study of the height, Sarnak conjectures that
it is a Morse function on the Teichm\"uller space (see \cite{Sa}). There is also a recent work by Sarnak and Str\"ombergsson regarding critical points of the height function on the space of ($n$-dimensional) flat tori \cite{SSt}.

 This paper relies heavily on the methods and techniques
developed by  Osgood, Phillips, and Sarnak \cite{OPS1} \cite{OPS2}
\cite{OPS3} (see also the good survey paper \cite{Sa2}  by Sarnak);
however, only the relatively easy part of their
analysis of Polyakov-Alvarez formula (see \eqref{Polyakov-D}) about the conformal effect of
the metric to its height is used.

\subsection*{Plan of the paper}
In Section~\ref{S:heat kernel}, we recall some basic results about heat kernels and heights.
Sections
\ref{Section hyperbolic},
\ref{S:separation}, \ref{S:proof collar}, and \ref{S:proof insertion} prove Theorem~\ref{proper on hyperbolic moduli intro}.
Sections
 \ref{S:conformal},
 \ref{S:uniform}, and \ref{S:continuity} explain and prove
Theorem~\ref{two uniformization homeo}.
Section~\ref{S:compactness} proves Theorem~\ref{compact isospectral}.
Section~\ref{S:flat proper} shows Theorem~\ref{proper on moduli II intro}.
Finally, we make some further remarks  in Section~\ref{S:remarks}.

\subsection*{Acknowledgment}
It is a great pleasure of the author to thank
Anton Deitmar, Leonid Friedlander,
Ezra Getzler, Chris Judge, Richard Melrose, Peter Sarnak, Scott
Wolpert, and Jared Wunsch for their remarks, suggestions, support,
and interest. He also thanks Ilia Binder, Bernard Gaveau, Victor
Ivrii, Carlos Kenig, Pei-hsun Ma, Robert McCann, Adrian Nachman, Jason Ng, and Jeremy
Quastel, for mathematical or nonmathematical helps.
Most part of this research was done
while the author was visiting Fields Institute at Toronto, IH\'ES at
Bure-sur-Yvette, and Max Planck Institute at Bonn; he is thankful for their hospitality.

The author acknowledges that this work is a result of his family's love, support, encouragement, and patience.
This paper is dedicated to them: his wife Dong-Soon Shim, and two daughters, Joyce Eunjin and Ashley Souyoung.

\section{Heat kernels and heights}\label{S:heat kernel}

In this section let $M$ be a compact Riemannian manifold
(possibly $\partial M \ne \emptyset$).
For a given Riemannian metric $\sigma$ on $M$, the Laplacian  is
the following second order elliptic operator on functions
(with Dirichlet boundary condition when $\partial M \ne \emptyset $):
in local coordinates,
\begin{align*}
  \Delta = \Delta_\sigma = - \frac{1}{\sqrt{\det(\sigma_{ij})}} \, \partial_i \,
  \sigma^{ij} \sqrt{\det(\sigma_{ij})} \  \partial_j   ,
\end{align*}
where $\sigma^{ij}$ is the inverse matrix of $\sigma_{ij}$.
We use $\partial^\sigma_n$, or just $\partial_n$, to denote
the outer normal derivative on the boundary.
 The key result of this
section is \eqref{double Dirichlet heat},
which is used to show the insertion lemma
(Proposition~\ref{hyperbolic insertion}).

\subsection{Heat kernels}

Let $P=P(x,y,t),\   t >0 $ be the heat kernel
(Dirichlet heat kernel when $\partial M \ne \emptyset$), i.e.
the fundamental solution of the heat equation
\begin{align*}
  & \partial_t P (x,y,t) + \Delta_x P(x,y,t) = 0 ,\\
  & \lim_{t \to 0 +} P(x,y,t) = \delta_x (y)  \\
  & (P(x, y , t) = 0 \ \ \hbox{for } x \in \partial M ).
\end{align*}
We use the convention that if $x,y$ belong to two different connected components, then $P(x,y)=0$.

One of the fundamental results of heat kernels is the estimate provided by
Cheeger, Gromov, and Taylor \cite{CGT} (see also \cite{Ch}).
\begin{theorem}\label{CGT estimate}
  For a complete Riemannian manifold $M^n$, and
  $x,y \in M$, $r >0$ such that the geodesic distance $ d(x,y) > 2r$,
  the following inequality holds.
  \begin{align*}
    P(x,y,t) \leq  c(n) (t^{-n/2} + t r^{-(n+2)})
    (\Phi(x, r) \Phi(y, r))^{-1/2}
    \exp(-\frac{(d(x,y)-2r)^2}{4t}).
  \end{align*}
  Here $\Phi(x,r)$ is the isoperimetric constant
  of the geodesic ball $B(x, r)$, i.e.
  \begin{align*}
    \Phi(x,r) = \inf_{\Omega} \frac{\vol(\partial
      \Omega)^n}{\vol(\Omega)^{n-1}} ,
  \end{align*}
  where $\Omega$ ranges over all open submanifolds
  which have compact closures with smooth boundary in $B(x,r)$.
\end{theorem}

\begin{remark}
Judge \cite{Ju} applies this heat kernel estimate to show
the convergence of heat kernels when the metric of a manifold degenerates.
Ji \cite{Ji2} uses a shaper estimate of Li and Yau \cite{LY}
under an additional assumption that the Ricci curvature is bounded from
below. Li and Yau's estimate fits into our situation; however, the weaker
but more general estimate of Cheeger, Gromov, and Taylor is
enough for our purpose.
\end{remark}
\begin{remark}\label{double}
    By a simple argument which resembles doubling,
    a similar estimate as given in Theorem~\ref{CGT estimate}  holds
    for the Dirichlet heat kernel of a compact Riemannian manifold with
    boundary.  Consider the following example.
    Let the closed cylinder $[0, 1] \times S^1$ have the
    metric $du^2 + h(u, v)dv^2 $, where $u$ is the parameter of $[0,1]$,
    $v$ is the parameter of $S^1$, and the function
    $h(u,v)$ is smooth and positive
    on the cylinder.  First extend the cylinder to a larger cylinder
    $[-1/2 , 3/2] \times S^1$ and the metric to $du^2 + H(u, v)dy^2$, where the
    function $H$ extends $h$ smoothly (or in such a way it has as much
    regularity as we need)  such that $H$ is constant near the new
    boundary $\{ -1/2\} \times S^1$ and $\{3/2\} \times S^1$. Then
    double the larger cylinder $[-1/2, 3/2] \times S^1$ to get a torus, say
    $T$. Note that the doubled metric on this torus is smooth.
    Consider the heat kernel $P^T$ of this torus $T$ with the newly
    constructed metric. Apply the heat kernel bound of
    Theorem~\ref{CGT estimate} to this heat kernel $P^T$. Now the original cylinder
    $[0,1] \times S^1$ is embedded in the torus $T$, and the Dirichlet heat
    kernel $P^D$ of the original cylinder $[0,1] \times S^1$ is bounded
    by the heat kernel $P^T$ of the torus by the maximum principle. Then we get
    \begin{align*}
        P^D (x,y, t) \le P^T (x,y, t) \le \hbox{const.}
        \exp( - \hbox{const.}/ t) \ \
        \ (0 < t < 1)
    \end{align*}
    where by Theorem~\ref{CGT estimate}
     the constants $(>0)$ depend only on the distance between $x$ and $y$ and
    on some appropriate metric balls about these two points.   The same inequality
    \begin{align}\label{double Dirichlet heat}
        P(x,y, t) \le \hbox{const.} \exp( - \hbox{const.} / t)  \ \
        \ (0 < t < 1)
    \end{align}
    with the same dependency of the constants as above,
    holds for a compact Riemannian manifold with piecewise smooth
    and Lipschitz boundary; in particular, it holds for a compact hyperbolic surface with piecewise geodesic boundary.
\end{remark}

Theorem~\ref{CGT estimate} and
Remark~\ref{double} are used later to give a proof of
the so-called insertion lemma (Proposition~\ref{hyperbolic insertion}).
\subsection{Heights}\label{Heights}
For a smooth compact Riemannian manifold $M$ with $\partial M \ne \emptyset$,
define the spectral zeta function
\begin{equation*}
  \zeta(s) = \frac{1}{\Gamma(s)} \int_0 ^\infty t^{s-1}  \int_M P(x,x,t)dx  dt,
  \ \ \ \Re s > \dim{M}/2 ,
\end{equation*}
where $dx$ is the Riemannian volume form. The integral $\int_M P(x,x,t) dx$ is called
the \emph{trace} of (Dirichlet) heat kernel, and we use $\Tr_D e^{-t\Delta}$ to denote it,
where the subscript $D$ means the Dirichlet condition.  
The spectral zeta function $\zeta (s)$ has its meromorphic extension to $\mathbb{C}$ and is holomorphic
at $s=0$. In two-dimension,
its regularity at $s=0$ can easily be seen
by using the asymptotic formula by McKean and Singer
(see \cite{OPS1} section 1, \cite{MS} (5.2))
\begin{align}\label{asymptotic formula}
    \int_M P(x,x,t) dx
    =& \frac{1}{4 \pi t} \text{Area} (M)
    -\frac{1}{8\sqrt{ \pi t}} \int_{\partial M} ds \\\nonumber
    & + \frac{1}{12 \pi} \int_M K dx
    - \frac{1}{12 \pi} \int_{\partial M} k ds
    + o(\sqrt{t}) \ \ \hbox{(as $t \to 0$)},
\end{align}
where $K$, $k$ denote the Gaussian curvature
and the boundary geodesic curvature, respectively,
and $ds$ denotes the arc length element.

Define the determinant $\det  \Delta$ of the Laplacain as following:
\begin{align*}
  -\log \det \Delta = \zeta ' (0).
\end{align*}
We call $-\log \det \Delta$ the \emph{height} of $M$ and denote it $h(M)$.

\begin{remark}
    It is not hard to see that the height $h$ gives
    a continuous function on the space of  Riemannian metrics
    in the $C^\infty$-topology.
    This is because for a given time $t >0$
    the heat kernel $P(x, y, t) $ depends
    continuously on the metric and
    \begin{align*}
        P(x, y, t) & = \int_M P(x,z,t/2) P(z,y,t/2) dz \\
        & \le e^{-\lambda_1 (t-T)} \sqrt{P(x,x,T)}
        \sqrt{P(y,y,T)}, \hbox{ \ \ for $t \ge T >0$},
    \end{align*}
    and
    $P(x, x, t) = O(t^{-\dim M /2})$
    as $t \to 0$,
    where the first eigenvalue $\lambda_1$ and the constant of $O(t^{-\dim M /2})$ depend
    continuously on the metric.
    Therefore,  $\zeta(s)$ and $\frac{d}{ds} \zeta(s)$
    for each $\Re s > \frac{\dim M}{2} $ are continuous
    with respect to the metric, and their analytic
    extensions at $s=0$ are continuous on the metric as well.
\end{remark}

\section{Heights of bordered hyperbolic surfaces}\label{Section hyperbolic}

Theorem~\ref{proper on hyperbolic moduli intro} is proven in this section.
Let $M$ be a hyperbolic surface (Gaussian curvature $-1$) with geodesic
boundary, i.e. $\partial M = \cup_{i=1}^n b_i $, where each $b_i$ is a closed
curve with zero geodesic curvature. Note that by Gauss-Bonnet theorem,
$\chi(M) < 0$. As typical in analysis of hyperbolic surfaces,
the so-called \emph{thick and thin decomposition} is used.
The insertion lemma
(Proposition~\ref{hyperbolic insertion}) is applied to this decomposition. The method of using the insertion lemma to study the height of hyperbolic surfaces is first considered by Sarnak \cite{Sa2}.

Let $\tilde{M}$ be the double of $M$, then
$\tilde{M}$ is a smooth closed hyperbolic surface. Let $g$ be the genus of $\tilde{M}$.
It is a well-known fact (see, for example, \cite{W})  that there is a
constant $1> c_* >0$, depending only on $g$, with the following property.
There are only finitely many (at most $3g-3$) closed
primitive geodesics of $\tilde{M}$, say,
$\gamma_1 , \cdots, \gamma_k$, whose
lengths $l(\gamma_i)$ are less than $c_*$ (these geodesics are called
\emph{short geodesics}); for each $\gamma_i$, there is a tubular
neighborhood $C_i$ called \emph{standard collar} \cite{W}, of width
\begin{align*}
  \sinh^{-1} ( 1/\sinh(\frac{1}{2} l (\gamma_i))) \approx  2
  \log \frac{2}{l (\gamma_i)};
\end{align*}
each $C_i$ is a hyperbolic cylinder
with the core geodesic $\gamma_i$, and these collars are all mutually disjoint.
The standard collar $C_\gamma$ of a short geodesic $\gamma$ of length $l$ can
be regarded as a domain
$[0,l] \times [l , \pi - l]$ with variables $(u, v)$ such that $u=0$ is
identified with $u=l$. The hyperbolic metric is $1/\sin^2(v) (du^2 + dv^2)$.
\emph{The standard subcollar} $SC_\gamma$ is defined
as the subset
\begin{align*}
\{ (u,v) \in C_\gamma \ | \ 2 l \le v \le \pi -2l \}/{(u=0 \sim u=l) } \ \ \  \subset C_\gamma.
\end{align*}
From the argument given by Wolpert
(see \cite{W}, section 2.6 and 2.7, especially p296),
the surface
$\tilde{M} \setminus \bigcup_{\gamma_i} SC_{\gamma_i}$ has
uniformly bounded geometry, i.e. the set
$\{ \tilde{M} \setminus \bigcup_{\gamma_i} SC_{\gamma_i} \}$
of such surfaces forms a compact
set in the $C^\infty$ topology of the space of Riemannian manifolds
modulo isometries. Moreover, there exists a constant $\delta^*
= \delta^*(c^*) > 0$ such that the tubular neighborhood
\begin{align*}
    A_\gamma=
    \{ x \in C_\gamma \ | \ \  \dist(x, \partial SC_\gamma )
    \le \delta^* \} \subset C_\gamma
\end{align*}
has uniformly bounded geometry.

For a short geodesic $\gamma$ in the double $\tilde{M}$, there are only three
possible cases:
\begin{align*}
  I: & \ \ \ \ \gamma \cap \partial M = \emptyset ,\\
  II: & \ \ \ \ \gamma \subset \partial M ,\\
  III: & \ \ \ \ \gamma \pitchfork \partial M.
\end{align*}
It is easy to see that in the second case,
$SC_\gamma \setminus \partial M$ consists of two isometric hyperbolic cylinders,
whose ends are $\gamma$ and a part of $\partial SC$. In the third case,
$SC_\gamma \setminus \partial M$ consists of two  isometric hyperbolic 4-gons,
each contains half of $\gamma$.

\begin{proposition}\label{heights of collars}Let
$SC_I$, $SC_{II}$, $SC_{III}$ denote the regions
$M \cap (SC_\gamma \setminus \partial M)$ corresponding to geodesics
$\gamma$ in cases $I, II, III$, respectively. Their heights satisfy:
  \begin{align}\label{SC I}
    h(SC_I) &\sim \frac{\pi^2}{6 l(\gamma) }
    + \log l(\gamma) + O(1) \ \ \ \ \hbox{\emph{by Lundelius \cite{Lu} and Sarnak \cite{Sa2}}},\\
    \label{SC II}
    h(SC_{II}) &\sim  \frac{\pi^2}{12 l(\gamma) }
    + \log l(\gamma) + O(1) , \\\label{SC III}
    h(SC_{III}) & \sim   \frac{\pi^2}{12 l(\gamma) }
    +O(1)  .
  \end{align}
\end{proposition}
\begin{proof}
The proof is given in Section~\ref{S:proof collar}.
\end{proof}

\begin{remark}
The domain $SC_{III}$ has right-angle corners; thus the previous definition of height in Section~\ref{Heights} does not directly apply. This
subtlety will be addressed in Section~\ref{S:separation}.
\end{remark}

\begin{proposition}\label{hyperbolic insertion}(\textbf{Insertion Lemma})
    Let $N$ and $A$ be the subsets $M \cap (\bigcup_{\gamma_i}
    SC_{\gamma_i})$ and $M \cap (\bigcup_{\gamma_i}
    A_{\gamma_i})$ of $\tilde{M}$, respectively.
    Then
    \begin{align}
        \label{define h}\hbox{the heights $h(N)$, $h(M\setminus N)$, and $h(A)$
        are defined;}\\\label{h decomp}
        \hbox{$h(M) \geq  h(N) + h(M \setminus N)+ O(1)$,}
    \end{align}
    where the constant depends only on $A$ and
    $\overline{\partial N \setminus \partial M}$.
\end{proposition}
\begin{proof}
    This result is a modified version of Osgood-Phillips-Sarnak's insertion lemma
    which was proved for the flat surface case  \cite{OPS3} (see also \cite{Kh}).
    The proof is given in Section~\ref{S:proof insertion}.
\end{proof}

$M\setminus N$ and $A$ have  uniformly bounded geometry and
so their heights are also bounded uniformly.
As a  consequence of
Proposition~\ref{hyperbolic insertion} and
Proposition~\ref{heights of collars},
we get the following asymptotic inequality.

\begin{theorem}\label{hyperbolic height} With the above setting,
  \begin{align*}
    h(M) \geq & \sum_{\gamma_i \cap \partial M = \emptyset} \big{(}
    \frac{\pi^2}{6l(\gamma_i)} + \log l(\gamma_i)\big{)}\\
    &  + \sum_{\gamma_i \subset \partial M}\big{(}\frac{\pi^2}{12 l(\gamma_i) }
    + \log l(\gamma_i)\big{)} \ \  + \sum_{\gamma_i \pitchfork \partial M}
    \frac{\pi^2}{12 l(\gamma_i) }
    \ \ +  O(1).
  \end{align*}
\end{theorem}
\begin{remark}
    An explicit formula of the height for geodesically bordered hyperbolic surfaces
    is given in terms of Selberg's zeta function by Bolte and Steiner
  \cite{BS}; however, it seems quite complicated to use their formula to get an
  asymptotic estimation. In the case of closed hyperbolic surfaces, Wolpert
  \cite{W} succeeded in using Selberg's zeta function expression of the height;
  his proof was later simplified by Lundelius \cite{Lu}, whose method is
  in the same spirit as ours in using the idea of insertion lemma. Bismut and Bost
  \cite{BB} took an alternative algebraic geometry approach.
\end{remark}

\begin{remark} One may try to refine the estimate in
  Theorem~\ref{hyperbolic height} by adding contributions from the low eigenvalues
  as in \cite{Lu} \cite{W}.
\end{remark}

\begin{corollary}\label{proper on hyperbolic moduli}
  (Theorem~\ref{proper on hyperbolic moduli intro})  On the moduli space
  $\mathcal{M}$ of compact hyperbolic
  (Gaussian curvature $\equiv -1$)  surfaces with
  geodesic boundary,   the  function $h$ is proper, i.e.
  \begin{align*}
    h(M) \rightarrow +\infty
  \end{align*}
  as the isometry class $[M]$ approaches $\partial \mathcal{M}$.
\end{corollary}
\noindent For a (finite dimensional) topological space $\mathcal{T}$, we say that a
sequence $t_i \to \partial \mathcal{T}$ if
the set $\{ t_i \}_{i=0} ^\infty$
is not contained in any compact subset of $\mathcal{T}$.

\section{Separation of variables, heat kernels, and heights}\label{S:separation}
Domains like $SC_{III}$ have points of special type singularity in their boundaries:
 there are four right-angle
corners in $\partial SC_{III}$. In this section we study the heat kernels and heights of
such domains. The subtlety due to their singularity
can be resolved by the separation of variables technique.
(See \cite{Ji1} for a different but more extensive use of separation of variables in studying the spectrum
of a Riemann surface.)

\subsection{$1$-dimensional heat traces}
Before we proceed the separation of variables, let's consider the $1$-dimensional case.

Let $[A, B] \subset \mathbb{R}$ be a finite closed interval with
metric $dx$. For a smooth function $\phi$ on $[A,B]$,
the Laplace operator of the metric $e^\phi dx$ is given by
\begin{align*}
  \Box_\phi = - e^{-2\phi} \big{[}\frac{d}{dx}\big{]}^2
  + \phi' e^{-2\phi} \frac{d}{dx}.
\end{align*}
Let $Q_\phi=
-e^{-2\phi} \big{[}\frac{d}{dx}\big{]}^2$, then
\begin{align}\label{Q phi}
Q_\phi = \Box_\phi -\phi'
e^{-2\phi} \frac{d}{dx}.
\end{align}

Let $e(x,y, t)$ be the fundamental solution of the heat equation of $Q_\phi$:
\begin{align*}
\frac{\partial}{\partial t} e(x,y,t) + [Q_{\phi}]_x e(x,y,t) = 0, \\
\lim_{t \to 0} e(x,y,t) = \delta_x (y),\\
e(x,y,t) = 0 \hbox{   for $x \in \{ A, B \}$},
\end{align*}
where $\delta_x$ is the Dirac $\delta$-function with respect to the metric $e^\phi dx$, 
i.e. $\int \delta_x (y) f(y) e^{\phi(y)}dy = f(x)$, for every smooth function $f$.

Let  $\varphi$ be an arbitrary smooth function on $[A,B]$. Define
\begin{align}\label{1-D trace}
\Tr_D \varphi e^{-tQ_\phi}=\int_A^B \varphi(x) e(x,x,t) e^{\phi(x)} dx ,
\end{align}
where the subscript $D$ means the Dirichlet boundary condition.
By using a result
  of McKean and Singer (\cite{MS}, pp 53. equation (5.2) and its proof, especially pp 55--56) 
  applied to \eqref{Q phi},
  we have
  \begin{align}\label{1-D asympt}
    \Tr_D \varphi e^{-t Q_\phi} = \frac{1}{\sqrt{4 \pi t}} \int_{A} ^B
   \varphi(x) e^{\phi(x)} dx
    - \frac{1}{4} ( \varphi(A) + \varphi(B))  + O(\sqrt{t})
    \ \ \hbox{(as $t \to 0$)}.
  \end{align}

Define the zeta function of $Q_\phi$:
\begin{align*}
  Z_\phi(s) = \frac{1}{\Gamma(s)} \int_0 ^\infty t^{s-1}
  \Tr_D e^{-t Q_\phi}dt.
\end{align*}
By \eqref{1-D asympt}, the regularization of $Z_\phi(s)$ at $s=0$ is enabled.
Now, by an argument identical to the derivation of the Polyakov formula
  (see for example, \cite{OPS1}, pp 155--156),
\begin{align}\label{1-D Polyakov}
  Z_\phi '(0) = - \frac{1}{2} (\phi(A) + \phi(B) )  + Z_0 ' (0).
\end{align}

\subsection{Separation of variables and heat kernels}

Let's first fix some notation.
\begin{definition}\label{D:C l [A,B]}
For $0 < l < \frac{\pi}{4}$ and $l \le  A < B \le \pi-l$, let $C^{l, [A,B]}$ be the smooth hyperbolic cylinder
defined as the domain
$[0,l] \times [A , B]$ with variables $(u, v)$, $u=0$
identified with $u=l$, and with the hyperbolic metric
$1/\sin^2(v) (du^2 + dv^2 )$. For example, $SC_I = C^{l, [2l , \pi -2l]}$,
$SC_{II}= C^{l, [2l , \pi/2]}$. Let $C_{III}^{l, [A,B]}$ be the subset
\begin{align*}
C_{III}^{l, [A,B]}= \{ (u, v) \in C^{l, [A,B]} \ |  0 \le u \le l/2 \} \ \
\subset C^{l, [A,B]}.
\end{align*}
For example, $SC_{III} = C_{III}^{l, [2l , \pi -2l]}$.
\end{definition}

We have the $L^2$-decomposition
\begin{align*}
  L^2 (C^{l,[A,B]} ) = L^2 ([A, B], \frac{1}{\sin^2(v)} dv)
  \oplus \big{[}\bigoplus_{m \in \mathbb{N}} L^2_{m, 1}(C^{l, [A,B]})
    \oplus  L^2_{m,2} (C^{l, [A,B]}) \big{]}
\end{align*}
where
\begin{align*}
  L^2_{m, 1}(C^{l, [A,B]}) & = \{f(v) \sqrt{\frac{2}{l}}\cos(2\pi m \,  u / l )\in L^2(C^{l, [A,B]}) \}
  \simeq L^2 ([A, B], \frac{1}{\sin^2(v)} dv),
  \\
  L^2_{m,2} (C^{l, [A,B]}) &=\{ f(v) \sqrt{\frac{2}{l}}\sin(2\pi m \,  u / l ) \in L^2(C^{l, [A,B]}) \}
  \simeq L^2 ([A, B], \frac{1}{\sin^2(v)} dv) .
\end{align*}

The Laplace operator
\begin{align*}
  \Delta
   =  - \sin ^2 (v) \Big{(} \frac{\partial^2}{{\partial u }^2}  +
   \frac{\partial^2}{{\partial v }^2}\Big{)}
\end{align*}
has  the decomposition
\begin{align*}
  (\Delta, L^2(C^{l, [A,B]}) )  & \simeq  \bigoplus_{m \in \mathbb{Z}}
  (\Delta_l (m), L^2 ([A, B], \frac{1}{\sin^2(v)} dv)),
\end{align*}
where
\begin{align*}
  \Delta_l (m)= -\sin^2 (v) \Big{(} \frac{\partial^2}{{\partial v }^2}
  - \frac{4 \pi^2 m^2}{l^2} \Big{)}.
\end{align*}
Denote  the Dirichlet boundary condition with subscript $D$, then
\begin{align*}
  (\Delta, L^2_D(C^{l, [A,B]}) )  & \simeq  \bigoplus_{m \in \mathbb{Z}}
  (\Delta_l (m), L^2_D ([A,B], \frac{1}{\sin^2(v)} dv)).
\end{align*}

Regarding $C^{l, [A,B]}_{III}$, it is easy to see that
\begin{align*}
    L^2 _D ( C_{III}^{l, [A,B]})
  & \simeq \bigoplus_{m \in \mathbb{N}}( L^2_{m, 2, D}(C^{l, [A,B]}),
  \frac{1}{\sin^2(v)} dv),\\
  (\Delta, L^2 _D ( C_{III}^{l, [A,B]})) &
  \simeq \bigoplus_{m \in \mathbb{N}}
  (\Delta_l (m), L^2_D ([A, B], \frac{1}{\sin^2(v)} dv)).
\end{align*}

Let $P^{l,[A, B]}(u_1, v_1 ; u_2, v_2 ; t)$, $P_{III}^{l,[A, B]}(u_1, v_1 ; u_2, v_2 ; t)$,
and $P^{l, m}(v_1, v_2; t )$ denote
the Dirichlet heat kernels of the operators
$ (\Delta, L_D^2(C^{l, [A,B]}) ) $, $ (\Delta, L_D^2(C_{III}^{l, [A,B]}) )$, and
$(\Delta_l (m), L^2_D ([A, B], \frac{1}{\sin^2(v)} dv))$, respectively.
We see that
\begin{align*}
& P^{l,[A, B]}(u_1, v_1 ; u_2, v_2 ; t)\\
& =
P^{l,0}(v_1, v_2;t)\frac{1}{l} \\
& +
\sum_{m \in \mathbb{N}}
P^{l, m}(v_1, v_2 ;t )
 \frac{2}{l} \Big{(}\cos(\frac{2\pi m u_1}{  l} ) \cos(\frac{2\pi m u_2 }{l})
+ \sin(\frac{2\pi m u_1}{ l} ) \sin(\frac{2\pi m u_2} {l} ) \Big{)}
\end{align*}
and
\begin{align*}
& P_{III}^{l, [A, B]}(u_1, v_1 ; u_2, v_2 ; t) =
\sum_{m \in \mathbb{N}}
P^{l, m} (v_1, v_2 ; t)
 \frac{4}{l} \sin(\frac{2\pi m u_1}{ l} ) \sin(\frac{2\pi m u_2} {l} ).
\end{align*}

\subsection{Heat traces, asymptotic expansions, and heights}\label{SS:asympt heat}
Let $\varphi(v)$ be a nonnegative smooth function on $[A,B]$ which is constant near the boundary points $v=A$ and $v=B$: for example, $\varphi \equiv 1$. Let's consider modified heat traces of the form
$\Tr \varphi e^{-t\Delta}$.
First,
\begin{align*}
 &\Tr_{L^2_D(C^{l, [A,B]})} \varphi e^{-t\Delta}\\
& = \int_A^B \int_0^l \varphi(v) P^{l, [A,B]} (u,v; u,v; t) \frac{1}{\sin^2(v)}du dv \\
& =  \int_A^B \varphi(v)P^{l,0}(v, v;t) \frac{1}{\sin^2(v)}dv
 +
2 \int_A^B  \varphi(v) \sum_{m \in \mathbb{N}}
P^{l, m}(v, v ;t )\frac{1}{\sin^2(v)}dv .
\end{align*}
Similarly,
\begin{align*}
\Tr_{L^2_D(C_{III}^{l, [A,B]})} \varphi e^{-t\Delta}
 = \int_A^B  \varphi(v) \sum_{m \in \mathbb{N}}
P^{l, m}(v, v ;t )\frac{1}{\sin^2(v)}dv.
\end{align*}
Thus, we see
\begin{align}\label{identity of heat kernels}
  2\Tr_{L^2 _D ( C_{III}^{l, [A,B]})} \varphi e^{-t \Delta}
  =     \Tr_{L^2_D(C^{l, [A,B]})}\varphi e^{-t\Delta}
  - \Tr_{L^2_D ([A, B])}\varphi e^{-t \Delta_l (0)},
\end{align}
where 
\begin{align*}
\Tr_{L^2_D ([A, B])}\varphi e^{-t \Delta_l (0)}
=\int_A^B \varphi (v) P^{l,0} (v,v;t) \frac{dv}{\sin^2(v)}.
\end{align*}

On the other hand, by the result (\cite{MS}  pp. 53--60, equation (5.2) and its proof)
of McKean and Singer,
it follows the asymptotic expansion
\begin{align}\label{asympt C}
&\Tr_{L^2_D(C^{l, [A,B]})}\varphi e^{-t\Delta}\\\nonumber
 &= \frac{l}{4 \pi t}  \int_A^B \varphi(v) \frac{dv}{\sin^2 (v)}
    -\frac{l}{8\sqrt{ \pi t}}  \Big{(}\frac{\varphi(A)}{\sin A} + \frac{\varphi(B)}{\sin B}
    \Big{)} \\\nonumber
    & + \frac{l}{12 \pi} \int_A^B  \frac{-\varphi(v)}{\sin^2 (v)} dv
    - \frac{l}{12 \pi}  ( \frac{\varphi(A) k (A) }{\sin A}
    + \frac{\varphi(B) k (B)}{\sin B})
    + o(\sqrt{t}) \ \ \hbox{(as $t \to 0$)},
\end{align}
where $k(A)$, $k(B)$ denote the constant geodesic curvatures of the boundaries $v=A$ and $v=B$
of $C^{l, [A,B]}$.

\begin{remark}\label{R:modified trace}
The above asymptotic expansion \eqref{asympt C} (up to $o(\sqrt{t})$)
uses that near the boundary the modifying function $\varphi$ is
constant along the normal direction. However,
the 1-dimensional expansion \eqref{1-D asympt} (up to $O(\sqrt{t})$) needs no such restriction on $\varphi$. We expect that for general case the asymptotic expansion would contain some derivatives of $\varphi$ (with respect to the normal direction to the boundary).
\end{remark}

For the asymptotic expansion of $\Tr_{L^2_D ([A, B])}\varphi e^{-t \Delta_l (0)}$,
apply \eqref{1-D asympt} with $\phi = -\log \sin v$ to get
\begin{align}\label{1-D asympt special}
\Tr_{L^2_D ([A, B])}\varphi e^{-t \Delta_l (0)}
 = \frac{1}{\sqrt{4 \pi t}} \int_{A} ^B
    \frac{\varphi(v)}{\sin v} dv
    - \frac{1}{4}(\varphi(A) + \varphi(B))  + O(\sqrt{t})
    \ \ \hbox{(as $t \to 0$)}.
\end{align}
\begin{remark}
Note that the fundamental solution $e$ in the integral \eqref{1-D trace} is now
$e(v,v,t)= P^{l,0} (v,v;t) \frac{1}{\sin v}$.
\end{remark}

From \eqref{identity of heat kernels}, \eqref{asympt C}, and \eqref{1-D asympt special},
an asymptotic expansion
\begin{align}\label{asympt C III}
\Tr_{L^2 _D ( C_{III}^{l, [A,B]})} \varphi e^{-t \Delta}=
\frac{C_1}{t} + \frac{C_2}{\sqrt{t}} + C_3 + O(\sqrt{t}) \ \ \hbox{(as $t \to 0$)}
\end{align}
follows;
this shows the spectral zeta function
\begin{align*}
Z_{III}(s)= \frac{1}{\Gamma(s)} \int_0 ^\infty t^{s-1}
\Tr_{L^2 _D ( C_{III}^{l, [A,B]})}  e^{-t \Delta} dt
\end{align*}
has regularization at $s=0$; thus, the well-definedness of the height $h(C_{III}^{l, [A,B]})$ follows.
By \eqref{identity of heat kernels},
\begin{align}\label{heights identity}
2 h(C_{III}^{l, [A,B]})& =
  h(C^{l, [A,B]})  - Z'(0).
\end{align}
Here,
\begin{align*}
  Z(s) = \frac{1}{\Gamma(s)} \int_0 ^\infty t^{s-1}
  \Tr_{L^2_D ([A, B])}e^{-t \Delta_l (0)}dt,
\end{align*}
and its regularization at $s=0$ is enabled by the asymptotic expansion \eqref{1-D asympt special}.

\section{Proof of Proposition~\ref{heights of collars}}\label{S:proof collar}
\eqref{SC I} was shown by Lundelius (\cite{Lu}, section 3.3) and Sarnak (\cite{Sa2}, Appendix)
by using the
Polyakov-Alvarez formula (see \eqref{Polyakov-D}) for a flat cylinder.  Exactly the same
proof as in \cite{Lu} can show \eqref{SC II}.

For \eqref{SC III}, we apply the results of Section~\ref{S:separation}.
Note that for a short geodesic of length $l$, the domains $SC_I$, $SC_{III}$ are
exactly $C^{l , [2l , \pi -2l]}$, $C_{III}^{l, [2l , \pi -2l]}$, respectively.
By \eqref{heights identity} for $A=2l , B=\pi -2l$,
\begin{align}\label{SC I and III} 2 \, h(SC_{III})& =
  h(SC_I)  - Z'(0).
\end{align}
By \eqref{1-D Polyakov} for $\phi = \log(1/\sin(v))$,
\begin{align*}
  Z' (0)= Z_\phi ' (0)  =\log( \sin(2l)) + Z_0'(0)
\end{align*}
where
$Z_0 '(0) = 2 \log\big{(}\frac{\pi- 4l}{\pi}\big{)} \zeta_R(0) + 2 \zeta_R ' (0)$
for the Riemann zeta function $\zeta_R (s)$.
Hence,
\begin{align*}
  Z'(0) \sim \log 2l + O(1),   \ \ \ \hbox{as } l \rightarrow 0  ,
\end{align*}
and so by \eqref{SC I and III} and \eqref{SC I} the proof of \eqref{SC III} is complete.

\section{Proof of Proposition~\ref{hyperbolic insertion} (Insertion Lemma)}\label{S:proof insertion}
The proof consists of two parts.
First, we establish that the definition of height
$h$ makes sense for $N$, $M\setminus N$, and $A$,
whose boundaries have corners of special type.
Then we show the inequality \eqref{h decomp}.

\subsection{Proof of \eqref{define h}}
Recall that $N=M \cap (\bigcup_{\gamma_i}
    SC_{\gamma_i})$ and $A$ is the tubular neighborhood $M \cap (\bigcup_{\gamma_i}
    A_{\gamma_i})$ of $\overline{\partial N \setminus \partial M}$. So,
the connected components of $N$ and $A$ are of the form
either $C^{l, [A,B]}$ or $C^{l, [A,B]}_{III}$
  (see Definition~\ref{D:C l [A,B]}).
The results of Section~\ref{S:separation} (especially \eqref{heights identity}) thus show that the values $h(N)$ and $h(A)$ are
 well-defined.

  We now show the well-definedness of the value $h(M \setminus N)$.
  Let $V$, $V_0$, $V_1$ and  $Q(x, y, t)$, $Q_0(x, y, t)$, $Q_1(x,y,t)$
  denote the domains $M\setminus N$, $M$, $(M\setminus N) \cap A$ and the corresponding Dirichlet heat kernels,
  respectively.
  It is enough to show $\int_V Q(x,x,t) dx $ has an asymptotic expansion
  of the form
  \begin{align}\label{asympt}
    \frac{C_1}{t} + \frac{C_2}{\sqrt{t}} + C_3 + O(\sqrt{t}) \ \ \hbox{(as $t \to 0$)}.
  \end{align}

  Choose a smooth
  function $\varphi$, $0 \leq \varphi \leq 1$ supported on $A$ as follows. On each connected component
  of $A$ (which is of the form either $C^{l, [A,B]}$ or $C^{l, [A,B]}_{III}$),
  $\varphi = \varphi(v)$ and $\varphi \equiv 1$ near $\overline{\partial N \setminus \partial M}$; $\varphi \equiv 0$ near the boundary $v=A$ and $v=B$.
Then
\begin{align*}
&\int_V Q(x,x, t) dx \\
&= \int_{V_1} \varphi Q_1(x,x,t) dx+ \int_{V_1} \varphi (Q(x,x,t) -Q_1(x,x,t)) dx  \\
& \ \ + \int_V (1-\varphi) (Q(x,x,t) -Q_0 (x,x,t)) dx +
\int_V (1-\varphi) Q_0 (x,x,t) dx\\
& = I_1 + I_2 + I_3 + I_4.
\end{align*}

To deal with $I_2$ and $I_3$, apply
 the maximum principle to $Q(x, \cdot, t) - Q_i (x, \cdot, t)$, $i=0,1$,
  and use \eqref{double Dirichlet heat} to see that
  \begin{align}\label{bound of heat difference}
    0 & \leq (1-\varphi(x)) (Q_0 (x,x,t)- Q (x,x,t)) \leq
    (1-\varphi(x)) c_1
    \exp(-\frac{c_2}{t}), \\\nonumber
    0 &\leq \varphi(x) (Q (x,x,t) - Q_1(x,x,t)) \leq \varphi(x) c_3
    \exp(-\frac{c_4}{t}), \ \  \ \ (0 <t<1).
  \end{align}
Thus, $\varphi(x)(Q(x,x,t)- Q_1 (x,x,t))$ and $(1-\varphi (x))(Q_0(x,x,t)- Q (x,x,t))$ decrease
  rapidly as $t \rightarrow 0$  uniformly for $x \in V$;
\begin{align*}
I_2 = o(t), \   I_3 = o(t) \ \ \hbox{as $t \to 0$}.
\end{align*}

Let's now consider $I_1$ and $I_4$. Note the connected components of $V_1$ are of the form
either $C^{l, [A,B]}$ or $C^{l, [A,B]}_{III}$. Apply \eqref{asympt C} and \eqref{asympt C III} to see $I_1$ has an asymptotic expansion of the form
\eqref{asympt}.
On the other hand, $I_4$ has an asymptotic expansion of the form \eqref{asympt} by the same reason as for \eqref{asympt C} (see Remark~\ref{R:modified trace}).
This completes the proof.

\subsection{Proof of \eqref{h decomp}}
  We follow the outline given in the proof of insertion lemma in \cite{OPS3}
  or in \cite{Kh}.
  Denote $M,\  N ,\  M\setminus N$ and the corresponding Dirichlet heat kernels
  by $M_0, M_1 , M_2$ and $P_0(x, y, t)$, $P_1(x, y, t), P_2(x,y,t)$,
  respectively.
  Pick a tubular neighborhood $A_1$ of
  $\overline{\partial M_1 \setminus \partial M_0}$
  with $A_1 \subset \subset A$, such that each connected component of  $A_1$ is of the form
  either $C^{l, [A,B]}$ or $C^{l, [A,B]}_{III}$ (see Definition~\ref{D:C l [A,B]}).
  Choose a compactly supported smooth
  function $\varphi$ on $A_1$, such that $0 \leq \varphi \leq 1$ and
  $\varphi \equiv 1$ near $\overline{\partial M_1 \setminus \partial M_0}$.
  Then
  \begin{align*}
    & h(M_0) - h(M_1) - h(M_2)\\
    & =  \frac{ d }{ds } \Big{|}_{s=0} \frac{1}{\Gamma(s)}\sum_{j=1}^2
    \int_0^1 \int_{M_j} t^{s-1} \varphi ( P_0(x,x, t) - P_j(x,x,t)) dx  dt
    \\  \ \ & \  \ + \frac{ d }{ds } \Big{|}_{s=0} \frac{1}{\Gamma(s)}
    \sum_{j=1}^2 \int_0 ^1  \int_{M_j} t^{s-1} (1- \varphi)(P_0(x,x, t) -
    P_j (x,x,t)) dx dt\\
    & \ \ + \frac{ d }{ds } \Big{|}_{s=0}
    \frac{1}{\Gamma(s)}
    \sum_{j=1}^2 \int_1 ^\infty  \int_{M_j} t^{s-1} \varphi(P_0(x,x, t) -
    P_j (x,x,t)) dx dt\\
    & \ \  + \frac{ d }{ds } \Big{|}_{s=0} \frac{1}{\Gamma(s)}\sum_{j=1}^2
    \int_1^\infty \int_{M_j} t^{s-1}(1- \varphi)
    ( P_0(x,x, t) - P_j(x,x,t) )dx  dt  \\
    & = (I) + (II)+ (III) + (IV)  .
  \end{align*}

 Let's consider $(I)$.  For the regions $A$, $A\cap M_1$, and $A \cap  (M_2)$, denote
  the corresponding Dirichlet heat kernels by
  $P^A _0 (x,y,t), P^A_1 (x,y,t)$, and $ P^A_2(x,y,t)$, respectively.
  Apply the maximum principle to $P_j (x, \cdot, t) - P_j ^A (x, \cdot, t)$
  and use \eqref{double Dirichlet heat}, to see that for $i=0,1,2$,
  \begin{align}\label{bound of heat difference}
    0 \leq \varphi(x) (P_j (x,x,t) - P_j ^A (x,x,t)) \leq \varphi(x) c_1(A)
    \exp(-\frac{c_2(A)}{t}), \ \  \ \ (0 <t<1).
  \end{align}
  Now, write
  \begin{align*}
    (I) =& \frac{ d }{ds } \Big{|}_{s=0}
    \frac{1}{\Gamma(s)}\sum_{j=1}^2 \int_0^1
    \int_{M_j} t^{s-1} \varphi(x) ( P_0^A(x,x, t) - P_j^A(x,x,t) dx  dt \\
    & +  \frac{ d }{ds } \Big{|}_{s=0}
    \frac{1}{\Gamma(s)}\sum_{j=1}^2   \int_0^1
    \int_{M_j} t^{s-1} \varphi(x) ( P_0(x,x,t) -P_0^A(x,x, t))dx dt \\
    &- \frac{ d }{ds } \Big{|}_{s=0}
    \frac{1}{\Gamma(s)}\sum_{j=1}^2 \int_0^1
    \int_{M_j} t^{s-1} \varphi(x)(P_j(x,x,t)- P_j^A(x,x,t) )]dx  dt  .
  \end{align*}
  By \eqref{bound of heat difference},
  the second and third terms have uniformly rapidly
  decreasing integrands supported inside $A$;
  the first term depends only on $A$ and
  $\overline{\partial M \setminus \partial N}$.
  Therefore we  see that
  \begin{align*}
    (I) =  O(1)
  \end{align*}
  with the constant dependent only on $A$ and
  $\overline{\partial M \setminus \partial N}$.

  It remains to show  $(II)$, $(III)$, and $(IV)$ are nonnegative.
  $(1-\varphi(x))(P_0(x,x,t)- P_j (x,x,t))$ decreases
  uniformly for $x \in M_j$ rapidly as $t \rightarrow 0$ (by
  the above argument as for $P_j -P_j ^A$) and $t \rightarrow \infty$,
  and so does $\varphi (x)(P_0(x,x,t)- P_j (x,x,t))$ as $t \rightarrow \infty$.
    Thus,
  \begin{align*}
    (II) &= \sum_{j=1}^2 \int_0^1  \int_{M_j} t^{-1}  (1-\varphi)
    (P_0(x,x, t) - P_j (x,x,t)) dx dt  \ \  \geq 0  ,\\
    (III) &= \sum_{j=1}^2 \int_1 ^\infty  \int_{M_j} t^{-1}
    \varphi(P_0(x,x, t) - P_j (x,x,t)) dx dt  \ \  \geq 0  ,\\
    (IV) &= \sum_{j=1}^2 \int_1 ^\infty  \int_{M_j} t^{-1}
    (1-\varphi)(P_0(x,x, t) - P_j (x,x,t)) dx dt  \ \  \geq 0 ,
  \end{align*}
  with the inequalities by the maximum principle.
  This completes the proof.

\section{Conformal metrics and heights}\label{S:conformal}
Let $\Sigma$ be a compact surface with boundary, and
$\sigma$, $\sigma_0$ be two metrics on $\Sigma$ conformal with each other, i.e.
$\sigma = e^{2\psi} \sigma_0$. Let $K, k$ (resp. $K_0, k_0$) be
the sectional curvature and the boundary geodesic curvature of
the metrics $\sigma$ (resp. $\sigma_0$). Then
\begin{align}\label{conformal curvature}
  K &=e^{-2\psi} (  \Delta_0 \psi    + K_0 ),\\
  \label{geodesic curvature} k &= e^{-\psi} ( k_0 + \partial_n \psi ),
\end{align}
where $\partial_n$ is the outer normal derivative for $\partial \Sigma$ with
respect to the metric $\sigma_0$.
We will use the subscript $0$ for the quantities of the metric $\sigma_0$.

\subsection{Polyakov-Alvarez formula}
Alvarez extended (see \cite{Al} pp 158--159 also see
\cite{OPS1}) the Polyakov formula to surfaces with boundary. Namely,
for the metric $\sigma_0$, the height $h$ satisfies
\begin{equation}\label{Polyakov-D}
  \begin{array}{ll}
    h(e^{2 \psi}\sigma_0 ) = &\frac{1}{6 \pi}
    \big{\{} \frac{1}{2}\int_\Sigma | \nabla_0 \psi |^2 \,dA_0 +
    \int_\Sigma K_0 \psi \, dA_0  + \int_{\partial \Sigma} k_0 \psi ds_0 \big{\}}\\ \\
    & + \frac{1}{4 \pi} \int_{\partial \Sigma} \partial_n \psi \,ds_0
    + h(\sigma_0)  .
  \end{array}
\end{equation}
In particular, for the scaling $\lambda^2 \sigma$, where $\lambda >0$,
\begin{align}\label{scale height}
  h(\lambda^2 \sigma) = \frac{\chi (\Sigma)}{3} \log
  \lambda     + h(\sigma).
\end{align}
\section{Uniform metrics}\label{S:uniform}
In this section let $\Sigma= \Sigma_{g, n}$ be a type $(g,n )$ compact orientable surface
with boundary, i.e. $\Sigma$  is obtained
by removing $n$ disjoint disks from a closed surface
of genus $g$. Assume  $\chi (\Sigma)= 2-2g-n < 0$.

\subsection{Uniform metrics of types I and II}
 The
following definitions are due to Osgood, Phillips, and Sarnak.
\begin{definition} \cite{OPS1}
  A metric $\sigma$ on $\Sigma$ is called a \emph{uniform metric of type I}
  if the resulting Riemannian manifold $M$ has constant
  sectional curvature and
  $\partial M$ is of zero geodesic curvature. The metric is called
  a \emph{uniform metric of type II} if  the resulting Riemannian manifold
  $M$ is flat, i.e. the sectional curvature is zero,
  and $\partial M$ is of constant geodesic curvature,
  with the same constant for all the boundary components.
\end{definition}
\begin{remark}
  For a uniform metric of type I, the area determines the sectional curvature of the surface;
  however, it does not determine, for a uniform metric of type II, the boundary
  geodesic curvature.
\end{remark}

\begin{definition}\label{space of uniform metrics}
  Let $ \mathbf{M}_I (\Sigma , A)$ (resp.
  $\mathbf{M}_{II}(\Sigma, A)$)
  denote the space of  uniform metrics of type I (resp. II) on
  $\Sigma$ with fixed area $A$.
  The topologies on these spaces are induced by the natural $C^\infty$-topology
  on the space of sections in the bundle
  $(T^* \Sigma)^2 $.
\end{definition}

\begin{definition}
  Define the moduli space $\mathcal{M}_I (\Sigma, A)$
  (resp. $\mathcal{M}_{II}(\Sigma, A)$ of metrics of type I
  (resp. II) as
  \begin{align*}
    \mathcal{M}_I (\Sigma, A) &= \mathbf{M}_I (\Sigma, A) /
    \Diff \,(\Sigma),  \\
    \mathcal{M}_{II} (\Sigma, A) &= \mathbf{M}_{II} (\Sigma, A) /
    \Diff \, (\Sigma),
  \end{align*}
  where   $\Diff \,(\Sigma)$ is the group of diffeomorphisms
  of $\Sigma$.
\end{definition}

\begin{remark}
  It is known that both $\mathcal{M}_I (\Sigma, A)$ and
  $\mathcal{M}_{II}(\Sigma, A)$ have real dimension $6g-6 + 3n$.
\end{remark}
\begin{remark}
  As discussed in \cite{OPS3}, a uniform metric of type II is
  nothing but a metric obtained from a closed flat surface with conical singularities,
  by removing
  the metric disks of a fixed radius centered at each conical points.
  It is known (\cite{OPS3}, Theorem 1.1)
  that the moduli space $\mathcal{M}_{II}(\Sigma, A)$ is in one-to-one
  correspondence with the subset $\mathcal{C}(\Sigma)'$ of
  the space $\mathcal{C}(\Sigma)$ of conical metrics----flat metrics
  with conical singularities----on $\Sigma$,
  such that for each metric in $\mathcal{C}(\Sigma)'$,
  the distances between conical points are $>2$. We refer
  the reader to \cite{OPS3} or \cite{Kh} for details of conical metrics.
\end{remark}

Let $\sigma$ be a uniform metric of type II with area $A$ on $\Sigma$.
By the uniformization theorem, we find a unique uniform metric
$\tau$ of type I with area $A$ within the conformal class of
$\sigma$ (see \cite{OPS1} \cite{Ab}):
\begin{align*}
  \sigma= e^{2 \psi} \tau
\end{align*}
where $\psi$ uniquely solves the normalized boundary value problem \cite{OPS1}
\begin{align*}\left\{%
  \begin{array}{ll}
    \Delta_\tau \psi  + \frac{2\pi \chi(\Sigma)}{A} =0, &
    \hbox{in $\Sigma$,} \\[1ex]
    \partial_n ^\tau \psi -
    \frac{2\pi \chi(\Sigma)e^{\psi}}{\int_{\partial \Sigma}
      e^\psi ds_\tau} =0  ,
    & \hbox{on $\partial \Sigma$,} \\[1ex]
    \int_\Sigma e^{2\psi} dA_\tau  = A.
  \end{array}%
  \right.
\end{align*}
This $\psi$ also solves uniquely
\begin{align*} \left\{%
  \begin{array}{ll}
    \Delta_\sigma \psi  + \frac{2\pi \chi(M) e^{-2\psi}}{A}  =0, &
    \hbox{in $\Sigma$,} \\[1ex]
    \partial_n^\sigma \psi -
    \frac{2\pi \chi(\Sigma)}{\int_{\partial \Sigma} ds_\sigma} =0  ,
    & \hbox{on $\partial \Sigma$,} \\[1ex]
    \int_\Sigma e^{-2\psi} dA_\sigma  = A .    \end{array}%
  \right. \end{align*}
Write this one-to-one correspondence between  uniform metrics of types I
and II as
\begin{align*}
  \Psi : \mathbf{M}_{II}(\Sigma, A) \rightarrow \mathbf{M}_I (\Sigma, A)  ,
  \ \Psi(\sigma) = \tau ,\\
  \Phi : \mathbf{M}_{I}(\Sigma, A) \rightarrow \mathbf{M}_{II} (\Sigma, A) ,
  \ \Phi(\tau) =\sigma .
\end{align*}
It induces a one-to-one correspondence
\begin{align*}
  \tilde{\Psi} : \mathcal{M}_{II}(\Sigma, A)
  \rightarrow \mathcal{M}_I (\Sigma, A)  , \ \  \tilde{\Phi} :
  \mathcal{M}_{I}(\Sigma, A) \rightarrow \mathcal{M}_{II} (\Sigma, A)  .
\end{align*}
This correspondence is nothing but a translation of  the uniformization
theorem of Osgood, Phillips, and Sarnak \cite{OPS1}  into our context.
We have the following theorem whose proof is given in the next section.
\begin{theorem}\label{Psi}
  The one-to-one maps $ \Psi, \Phi$ are continuous and therefore
  homeomorphisms.
\end{theorem}
As a direct corollary, we have:
\begin{theorem}\label{T:tilde Psi}
  The maps $ \tilde{\Psi}, \tilde{\Phi}$ are homeomorphisms.
\end{theorem}

\subsection{Properness of the height on $\mathcal{M}_I(\Sigma, A)$}
We can rewrite Corollary~\ref{proper on hyperbolic moduli} as the following.
\begin{theorem}\label{T:proper on moduli I}
  For each $A >0$,
  the height $h$ is a proper function on the moduli space
  $ \mathcal{M}_I (\Sigma, A)$, i.e.
  \begin{align*}
    h(M) \rightarrow + \infty
  \end{align*}
  as the isometry class $[M]$ approaches $\partial \mathcal{M}_I (\Sigma, A)$.
\end{theorem}

\section{Proof of Theorem~\ref{Psi}}\label{S:continuity}
Throughout this section we extensively use the proof of
the uniformization theorem of Osgood, Phillips, and Sarnak
(see especially pp.158--163 of \cite{OPS1}, but note that their
Laplacian and our Laplacian have opposite signs).
Let $\tau \in \mathbf{M}_I (\Sigma, A)$,
$\sigma \in \mathbf{M}_{II} (\Sigma, A)$. Suppose $\Psi (\sigma) = \tau$ (and so
$\Phi (\tau) = \sigma$).
We assume without loss of generality that
the fixed area $A$ is $-2\pi \chi (\Sigma)$ and so the constant
curvature $K_\tau$ is $-1$.

\subsection{Continuity of $\Psi$}
First normalize $\sigma$ by a conformal factor $e^{2\varphi_0}$,
where $\varphi_0$ is a  solution
of the following differential equation:
\begin{align*}
  \partial_n ^\sigma \varphi_0 + k_\sigma & = 0  \ \ \text{on} \
  \partial \Sigma  .
\end{align*}
The resulting metric $\sigma_0 = e^{2\varphi_0}\sigma $ has
vanishing geodesic curvature of the boundary, i.e. $k_{\sigma_0}=0$. Moreover,
it is possible to pick $\varphi_0$ so that it (and also $\sigma_0$) depends
continuously on $\sigma$ in $C^\infty$.
Construct $\varphi_0$ in the following way.
Let $\eta$ be a smooth function compactly supported on
the interval $[0, 1)$ and $\eta \equiv 1$ on $[0, 1/2]$.
For $s >0$, define  $H_s : \partial \Sigma \rightarrow \Sigma$
by
\begin{align*}
  H_s (p) = \exp_p ( s \vec{n} ), \  p \in \partial \Sigma ,
\end{align*}
where $\exp$ and $\vec{n}$ are the exponential map and the inward
normal vector with respect to the metric $\sigma$, respectively.  Let
\begin{align*}
  r:=\sup \{ s >0 \ |  \hbox{$H_t$ is an embedding for every $t < s$} \}.
\end{align*}
On each inward normal geodesic ray emanating from $p \in \partial \Sigma$
of length $r$, first solve  the ODE
\begin{align*}
  f_p \, ' (t) + k_\sigma =0, \  f_p (0) =0  .
\end{align*}
Define $\varphi_0$ on $\Sigma$  as
\begin{align*}
  \varphi_0 (x) = \left\{%
  \begin{array}{ll}
    \eta( s/r ) f_p (s)
    , & \hbox{if $x= H_s (p)$ for some $s < r $ and
      $p \in \partial \Sigma$;} \\
    0, & \hbox{otherwise,} \\
  \end{array}%
  \right.
\end{align*}
then it has the desired property.

Without loss of generality, further normalize $\sigma_0$ so
that it has area $A$.

Consider the following functional with the reference metric $\sigma_0$:
\begin{align*}
  F_1 (\varphi) & = \frac{1}{2} \int_\Sigma | \nabla_0 \varphi |^2 dA_0
  + \int_{ \Sigma} K_0 \varphi dA_0
  -\pi \chi(\Sigma) \log \big{(}\int_{ \Sigma} e^{2\varphi} dA_0 \big{)}.
\end{align*}
Subject to the constraint
\begin{align}\label{constraint 1}
  \int_\Sigma \varphi dA_0 &= 0,
\end{align}
 $F_1$   is strictly convex and there exists a unique minimizer
$\psi$ of $F_1$ such that after a proper rescaling,
the metric $\tau = e^{2 \psi} \sigma_0 = \Psi(\sigma)$ is in $\M_{I}(\Sigma, A)$ (see \cite{OPS1}).
By the first variation of $F_1$, we see that $\psi$ satisfies
\begin{align}\label{PDE of psi}
  \Delta_0 \psi & =- K_0 + \frac{2\pi \chi (M)}{\int_M e^{2\psi}dA_0} \,
  e^{2\psi}  \ \ \hbox{in $\Sigma$},    \\
  \label{Neumann psi}
  \partial_n \psi  & = 0  \ \ \hbox{on $\partial \Sigma$}
\end{align}
(see \cite{OPS1} pp161).
For the continuity of $\Psi$, we show next that
$\psi$ depends continuously on $\sigma$ in $C^\infty$.

Suppose there is a sequence of metrics $\{\sigma ^i \}$
such that $\sigma^i \rightarrow \sigma$
in $C^\infty$. Consider the corresponding normalized metrics
$\sigma_0^i$ constructed as above, and see that
$\sigma_0^i \rightarrow \sigma_0$ in $C^\infty$.
Let $M_i$ (resp. $M$) denote
the Riemannian manifold $(\Sigma, \sigma_0^i) $
(resp. $(\Sigma, \sigma_0) $), and let
$F_1 ^i$ be the corresponding functional as above,
and $\psi ^i$  be the corresponding unique minimum, subject
to the constraint \eqref{constraint 1} with respect to $\sigma_0 ^i$.

We now show that $\psi^i \to \psi$ in $C^\infty$. Several
steps are similar (sometimes verbatim) to the arguments
given in \cite{OPS1}. In the following steps, if not specified, the metric
related quantities such as measure or gradient inside the integral
$\int_{M_i}$ (resp. $\int_{M}$) will be those of the corresponding metric
$\sigma_0 ^i$ (resp. $\sigma_0$); it is also important  to keep in mind that all these
metric related quantities are comparable with each other, respectively,
  because $\sigma_0^i \to \sigma_0$
in $C^\infty$.

\subsubsection*{Step 1.}
Since $\sigma_0 ^i \rightarrow \sigma_0 $ in $C^\infty$, it is clear that
\begin{align}\label{F i bound}
  F_1 ^i ( \psi^i ) \leq F_1 ^i (0 ) \leq \text{const.}.
\end{align}

\subsubsection*{Step 2: Obtain a priori bounds
of $\psi^i$.}
By Poincar\'e inequality  (with the constraint \eqref{constraint 1}),
\begin{align*}
  \Big{|}\int_{M_i} K_0 ^i \psi^i \Big{|}  \lesssim
  \Big{(}{\int_{M_i} | \psi^i |^2  }\Big{)}^{1/2}   \lesssim
  \Big{(} \int_{M_i} | \nabla \psi ^i |^2 \Big{)}^{1/2}.
\end{align*}
By Jensen's inequality and noting $\int_{M_i} \psi^i = 0$,
\begin{align*}
  \log \big{(}{\int_{M_i} \exp({2 \psi^i }})\big{)} \geq \log A.
\end{align*}
Combining these inequalities with \eqref{F i bound}, we get
\begin{align*}
  \int_{M_i } |\nabla \psi^i |^2   \leq \text{const.}.
\end{align*}
Since $\sigma_0 ^i \rightarrow \sigma_0 $ in $C^\infty$, we also have
\begin{align}\label{grad bound for Psi}
  \int_{M } |\nabla \psi^i |^2   \leq \text{const.}.
\end{align}

\subsubsection*{Step 3:
$\psi^i \to \psi $ in the space $L^2 (M)$.}
\eqref{grad bound for Psi} implies
\begin{align*}
\int_M | \psi^i |^2 \lesssim \int_{M} | \nabla \psi ^i |^2 \leq \hbox{const.}.
\end{align*}
By Rellich's compactness,  there exists $\psi^\infty$ in
the Sobolev space $W^1 (M)$ such that a subsequence
$\psi^{i_k} \rightarrow \psi^\infty $ weakly in $W^1 (M)$, strongly
in $L^2 (M)$, and pointwise a.e.
Therefore
\begin{align*}
  \int_M |\nabla \psi^\infty |^2  & \leq \liminf \int_{M_{i_k}}
  |\nabla \psi^{i_k}|^2 , \\
  \int_M K_0 \psi^\infty & = \lim \int_{M_{i_k}} K_0 ^{i_k} \psi^{i_k} ,\\
  \int_M \exp({2\psi^\infty}) & \leq \liminf \int_{M_{i_k}}
  \exp({2 \psi^{i_k} }),
\end{align*}
hence
\begin{align}\label{F 1 liminf}
  F_1 (\psi^\infty ) \leq \liminf F_1 ^{i_k} (\psi ^{i_k}) .
\end{align}
Let $\epsilon_k = \int_{M_{i_k}} \psi$,
$\psi_{\epsilon_k} = \psi- \frac{\epsilon_k}{A }$, therefore
 $\int_{M_{i_k}}
\psi_{\epsilon_k} = 0 $ . Obviously
$ F_1 ^{i_k} (\psi^{i_k}) \le F_1 ^{i_k} (\psi_{\epsilon_k}) $, and
$ F_1 ^{i_k} (\psi_{\epsilon_k}) \to F_1 (\psi)$.
Then by
\eqref{F 1 liminf} we have
\begin{align*}
  F_1 (\psi^\infty) \le    F_1 (\psi).
\end{align*}
This implies $\psi^\infty = \psi$ by the uniqueness of
the $F_1$ minimum.
Since this convergence holds for every convergent subsequence,
we conclude
that
\begin{align}\label{W 1 convergence}
  \psi^i \to \psi \ \ \hbox{ in $L^2 (M)$}.
\end{align}

\subsubsection*{Step 4.} For convergence in higher Sobolev spaces,
the elliptic regularity of
\eqref{PDE of psi}\eqref{Neumann psi} is used.
First, there is an \emph{a priori} bound
\begin{align}\label{a priori Sobolev in M}
  \| \psi^i \|_t \leq \hbox{const.}_t
\end{align}
for each ($t >0$) Sobolev $t$-norm  $\| \cdot \|_t$ on $M$.
(Trudinger's inequality is used when we bound the righthand-side
of \eqref{PDE of psi} for $\psi^i$ in $\| \cdot \|_0$. See \cite{OPS1}: equations
(3.19) and (3.24), and the remark below (3.29) in it.)
By Rellich's compactness and a diagonal argument, we see from
\eqref{W 1 convergence}, \eqref{a priori Sobolev in M} that
\begin{align*}
  \| \psi^i - \psi \|_t \to 0  \ \ \hbox{for each $t >0$}.
\end{align*}
This completes the proof of the $C^\infty$ convergence of
$\psi^i$ to $\psi$, and so the continuity of $\Psi$.

\subsection{Continuity of $\Phi$}
A similar
(sometimes verbatim) argument as in the case of $\Psi$ is used.

First find a solution for
\begin{align*}
  \Delta_\tau \varphi_0 = 1 \ \text{in $\Sigma$},
\end{align*}
so  the resulting metric
$\tau_0 = e^{\varphi_0}\tau$ is flat, i.e.
$K_{\tau_0} = 0$,
and depends continuously on $\tau$ in $C^\infty$.
$\varphi_0$ can be constructed in the following way.
Double
$\Sigma$ with metric $\tau$ to get a closed hyperbolic
surface $\tilde{\Sigma}$. $\Sigma$ is one half
of $\tilde{\Sigma}$. Let $\eta$ be a smooth
function compactly supported on the interval $[0, 1)$
and $\eta \equiv 1$ on $[0, 1/2]$.
For $s >0$, define  $H_s : \partial \Sigma \rightarrow \tilde{\Sigma}$
by
\begin{align*}
  H_s (p) = \exp_p ( s \vec{n} ), \  p \in \partial \Sigma ,
\end{align*}
where $\exp$, $\vec{n}$ are the exponential map and the outward
normal vector (for $\partial \Sigma$) with respect to $\tau$, respectively. Let
\begin{align*}
  r:=\sup \{ s >0 \ |  \hbox{$H_t$ is an embedding for every $t < s$} \}.
\end{align*}
Define $f_0$ on $\tilde \Sigma$ as
\begin{align*}
  f_0 (x) = \left\{%
  \begin{array}{ll}
    1, & \hbox{if $x \in \Sigma$;}\\
    \eta(  s/r )                , & \hbox{if $x= H_s (p)$ for some $s < r $
      and $p \in \partial \Sigma$;} \\
    0, & \hbox{otherwise.} \\
  \end{array}%
  \right.
\end{align*}
Now fix $\delta >0$ and a point $p$ outside the $r/2$-neighborhood
$U$ of $\Sigma$ so that
the $\delta$-geodesic ball of $p$ does not intersect $U$.
Pick a radially symmetric $C^\infty$ bump function,
say $f_1$, supported on this $\delta$-geodesic ball so that
the function $f= f_0 -  f_1$ satisfies
$\int_\Sigma f dA_\tau = 0 $. Then find the solution
for $\tilde{\varphi_0}$ (uniquely up to constant) where
\begin{align*}
  \Delta_\tau \tilde{\varphi_0} = f  \ \hbox{on $\tilde{\Sigma}$}.
\end{align*}
By restricting $\tilde{\varphi_0}$ to $\Sigma$ we find
the desired function $\varphi_0$.

From this point on, without loss of generality, further rescale
the flat metric $\tau_0$ so that it has area $A$.

Consider the following functional with the reference metric $\tau_0$.
\begin{align*}
  F_2 (\varphi) = \frac{1}{2} \int_\Sigma | \nabla_0 \varphi |^2 dA_0
  + \int_{\partial \Sigma} k_0 \varphi ds_0
  -2 \pi \chi(\Sigma) \log \big{(}\int_{\partial \Sigma}
  e^{\varphi} ds_0 \big{)},
\end{align*}
Subject to the constraint
\begin{align}\label{constraint 2}
  \int_{\partial \Sigma}  \varphi ds_0 = 0,
\end{align}
$F_2$ is strictly convex and there exists a unique minimizer
$\phi$ of $F_2$ so that the metric $e^{2 \phi} \tau_0$ is a
uniform metric of type II \cite{OPS1}. In particular, $\phi$ is
a harmonic function with respect to the flat  metric $\tau_0$.
We obtain $\sigma=\Phi(\tau)$ by rescaling this metric
(i.e. adding a constant to $\phi$) so that it has
the area $A$.
Let $T$ be the Dirichlet-Neumann operator of the metric $\tau_0$,
i.e. $T$ is the linear operator on functions on $\partial \Sigma$ given by
\begin{align*}
  T \varphi = \partial_n \tilde{\varphi}
\end{align*}
where $\tilde{\varphi}$ is the harmonic extension of $\varphi$ into
$\Sigma$.  As in \cite{OPS1}, $T$ is a positive self-adjoint
pseudodifferential operator of order $1$ on the space of
functions $\varphi$ with \ref{constraint 2}, and
we have
\begin{align}\label{norm of T}
  | T^{1/2} \varphi |_0 \sim  | \varphi |_{1/2}
  \ \ \hbox{and} \ \ |T \varphi |_0 \sim | \varphi|_1 ,
\end{align}
where $| \cdot |_t$ is the Sobolev $t$-norm on $\partial \Sigma$
with respect to the metric $\tau_0$.
Now as in \cite{OPS1} the minimizer $\phi$ of $F_2$ satisfies
\begin{align}\label{phi T equation}
  T \phi = -k_0 + \frac{2\pi \chi(\Sigma)e^\phi}{\int_{\partial \Sigma}
    e^\phi ds_0 },
\end{align}
and from harmonicity of $\phi$,
\begin{align}
  F_2(\phi) & = \frac{1}{2} \int_{\partial \Sigma} \phi T \phi  ds_0 +
  \int_{\partial \Sigma} k_0 \phi ds_0
  -2 \pi \chi(\Sigma) \log \big{(}\int_{\partial \Sigma}
  e^{\phi} ds_0 \big{)},\\ \label{F 2 boundary}
  &  =  \frac{1}{2} |T^{1/2} \phi |_0^2  +
  \int_{\partial \Sigma} k_0 \phi ds_0
  -2 \pi \chi(\Sigma) \log \big{(}\int_{\partial \Sigma} e^{\phi} ds_0 \big{)}.
\end{align}
By the expression of \eqref{F 2 boundary}, the functional $F_2$ induces a functional
on the space of functions on $\partial \Sigma$.
$\phi \big{|}_{\partial \Sigma}$ is the unique minimizer of
this induced functional subject to the constraint \eqref{constraint 2}.

For the continuity of $\Phi$, we show that
$\phi$ depends continuously on $\tau$ in $C^\infty$.
Our plan is first to derive this continuous dependence for
the restriction $\phi\big{|}_{\partial \Sigma}$ to the boundary
and then use harmonicity of $\phi$ to obtain the continuous dependence for
the whole function $\phi$.

Suppose there is a sequence of metrics $\{\tau ^i \}$
such that $\tau ^i \rightarrow \tau$
in $C^\infty$. Consider the normalized flat metrics
$\tau_0^i$, constructed as above, and  see that
$\tau_0^i \rightarrow \tau_0$. Let
$M_i$ (resp. $M$) be the Riemannian manifold
$(\Sigma, \tau_0^i) $ (resp. $(\Sigma, \tau_0)$).
Let $F_2 ^i$ be
the corresponding functional as above; and let $\phi ^i$
be the corresponding unique minimum subject to the constraint
\eqref{constraint 2} with respect to $\tau_0 ^i$.
To show that $\phi^i \to \phi$ in $C^\infty$, several
steps which are similar (sometimes verbatim) to the arguments
given in \cite{OPS1} are taken. In the following, if not specified, the metric
related quantities such as measure, gradient, or $T$ inside the integrals
$\int_{M_i}$,$\int_{\partial M_i}$ (resp. $\int_{M}$, $\int_{\partial M}$)
will be those of the corresponding metric
$\tau_0 ^i$ (resp. $\tau_0$); it is  also important to keep in mind that
all these metric related quantities are all comparable with each other, respectively, because $\tau_0^i \to \tau_0$ in
$C^\infty$.

\subsubsection*{Step 1.}
Since $\tau_0 ^i \rightarrow \tau_0 $ in $C^\infty$, it is clear that
\begin{align}\label{F 2 i bound}
  F_2 ^i ( \phi^i ) \leq F_2 ^i (0) \leq \text{const.} .
\end{align}

\subsubsection*{Step 2.} Use Step 1, to obtain some \emph{a priori} bounds of
$\phi^i$.
By \eqref{norm of T},
\begin{align*}
  \Big{|}\int_{\partial M_i} k_0 ^i \phi^i \Big{|}  \lesssim
  \Big{(}{\int_{\partial M_i} | \phi^i |^2  }\Big{)}^{1/2}   \lesssim
  \Big{(} \int_{\partial M_i} | T^{1/2} \phi ^i |^2 \Big{)}^{1/2}.
\end{align*}
By Jensen's inequality and noting $\int_{\partial M_i} \phi^i = 0$,
\begin{align*}
  \log \big{(}{\int_{\partial M_i} \exp({ \phi^i }})\big{)} \geq
  \log \int_{\partial M_i} 1 \ge \text{const.}.
\end{align*}
Combine these inequalities with \eqref{F 2 boundary} and \eqref{F 2 i bound}, we get
\begin{align*}
  \int_{\partial M_i } |T^{1/2} \phi^i |^2   \leq \text{const.}.
\end{align*}
Since $\tau_0 ^i \rightarrow \tau_0 $ in $C^\infty$, we also have
\begin{align}\label{T bound for Phi}
  \int_{\partial M } |T^{1/2} \phi^i |^2   \leq \text{const.}.
\end{align}

\subsubsection*{Step 3.} In this step, we show
$\phi^i \to \phi $ in the space $L^{2} (\partial M)$.
\eqref{T bound for Phi} combined with \eqref{norm of T} implies
\begin{align*}
\int_{\partial M} |
\phi^i |^2 \lesssim \int_{\partial M}  | T^{1/2}\phi ^i |^2 \leq \hbox{const.} .
\end{align*}
 By Rellich's compactness there is $\phi^\infty$
in  the Sobolev space $W^{1/2} (\partial M)$ such that a subsequence
$\phi^{i_k} \rightarrow \phi^\infty $ weakly in $W^{1/2} (\partial M)$ but
strongly in $L^2 (\partial M)$, and pointwise a.e.
This results
\begin{align*}
  \int_{\partial M} |T^{1/2} \phi^\infty |^2
  & \leq \liminf \int_{\partial M_{i_k}} |T^{1/2} \phi^{i_k}|^2 , \\
  \int_{\partial M} k_0 \phi^\infty & =
  \lim \int_{\partial M_{i_k}} k_0 ^{i_k} \phi^{i_k} ,\\
  \int_{\partial M} \exp({\phi^\infty})
  & \leq \liminf \int_{\partial M_{i_k}}
  \exp({ \phi^{i_k} }),
\end{align*}
therefore
\begin{align}\label{F 2 liminf}
  F_2 (\phi^\infty ) \leq \liminf F_2 ^{i_k} (\phi ^{i_k}).
\end{align}
Denote $\epsilon_k = \int_{\partial M_{i_k}} \phi$ and
let
$\phi_{\epsilon_k} = \phi- \epsilon_k /(\int_{\partial M_{i_k}} 1 )$
so that $\int_{\partial M_{i_k}} \phi_{\epsilon_k} = 0 $.
Then obviously
$ F_2 ^{i_k}( \phi^{i_k}) \le F_2 ^{i_k} (\phi_{\epsilon_k}) $ and
$ F_2 ^{i_k} (\phi_{\epsilon_k}) \to F_2 (\phi)$.
From \eqref{F 2 liminf},
\begin{align*}
  F_2 (\phi^\infty) \le    F_2 (\phi).
\end{align*}
This implies $\phi^\infty = \phi$ by the uniqueness of the minimum
of $F_2$.
Since this convergence is for every convergent subsequence,
we conclude that
\begin{align}\label{W 1/2 convergence phi}
  \phi^i \to \phi \ \ \hbox{ in $L^{2} (\partial M)$}.
\end{align}

\subsubsection*{Step 4.} For convergence in higher Sobolev spaces
$W^{t} (\partial M), \ t > 0$, we use the elliptic regularity of
\eqref{phi T equation}.
We need the following lemma.
\begin{lemma}\emph{(See \cite{OPS3}: eqn. (2.10) and Lemma 2.5.  See also \cite{OPS1}: Lemma 3.5 and
Lemma 3.9.) }
  For a fixed flat metric on $\Sigma$, there exist positive constants  such that
  for any smooth function $f$ on $\Sigma$,
  \begin{align}\label{1/2 sobolev}
    |f|_{1/2}^2  \leq \text{\emph{const.}}
    \big{\{} \int_\Sigma |\nabla f |^2 dA +
    \big{(}\int_{\partial \Sigma} f ds \big{)}^2 \big{\}} ,
    \\
    \label{Trudinger inequality}
    \int_{\partial \Sigma} e^f \frac{ds}{\int_{\partial \Sigma} ds} \leq
    \text{\emph{const.}} \  \exp \Big{(} \text{\emph{const.}}
    \int_\Sigma |\nabla f |^2 dA + \frac{\int_{\partial \Sigma} f ds}{\int_{\partial \Sigma} ds} \Big{)}.
  \end{align}
    In particular,
  \begin{align} \label{Trudinger T bound}
    \int_{\partial \Sigma} e^f  \frac{ds}{\int_{\partial \Sigma} ds} \leq \hbox{\emph{const.}}
      \
    \exp \Big{(} \hbox{\emph{const.}} | T^{1/2} f |_0^2 +
    \frac{\int_{\partial \Sigma} f ds }{\int_{\partial \Sigma} ds} \Big{)} .
  \end{align}
\end{lemma}

First, it follows that
\begin{align}\label{a priori Sobolev phi on partial M}
  | \phi^i |_t \leq \hbox{const.}_t  \ \
   \hbox{for each $ t> 0$}.
\end{align}
We have used \eqref{T bound for Phi} and \eqref{Trudinger T bound}
when we bound
the righthand-side of \eqref{phi T equation} for $\phi^i$ in $| \cdot |_0$.

As in the Step 4 for the case of $\Psi$,
apply Rellich's compactness and a diagonal argument to
\eqref{W 1/2 convergence phi} and
\eqref{a priori Sobolev phi on partial M} to get
\begin{align}\label{boundary convergence}
  | \phi^i - \phi |_t \to 0  \ \ \hbox{for each $t >0$}.
\end{align}

\subsubsection*{Step 5}
By elliptic regularity theory, observe the following standard fact.
\begin{lemma}\label{elliptic harmonic}
  If
  \begin{align*}
    \Delta f = 0 \text{   in $\Sigma$}\ \ \text{ and  }  \ \
    |f|_t <\infty \text{  on $\partial \Sigma$},
  \end{align*}
  then
  \begin{align*}
    \| f \|_t \leq \text{\emph{const.}}_t | f |_{t -1/2}
    \leq \text{\emph{const.}}_t |f|_t .
  \end{align*}
\end{lemma}
Apply this lemma to $\phi^i$ to get
$\| \phi^i \|_t \le \hbox{const.}_t$. By Rellich's compactness
and a diagonal argument, there is a subsequence
$\phi^{i_k}$ converging to some function $\phi^\infty$ on
$\Sigma$ in $C^\infty$.
This $\phi^\infty$ is harmonic for the metric $\tau_0$;
and by \eqref{boundary convergence},
$\phi^\infty \big{|}_{\partial \Sigma} = \phi\big{|}_{\partial \Sigma}$.
Therefore by uniqueness of Dirichlet problem we see
$\phi^\infty = \phi$ in $\Sigma$. This is for every convergent subsequence
and so we conclude that $\phi^i \to \phi$ on $\Sigma$
in $C^\infty$. The proof of
the continuity of $\Phi$ is thus complete.

\section{Compactness of the set of isospectral flat surfaces
with boundary}\label{S:compactness}
In this section we prove the following theorem.
\begin{theorem}\label{T:isospectral compactness}
  Let $\Sigma$ be a compact orientable surface with boundary. Assume $\chi(\Sigma) < 0$.
  Then
  each Dirichlet isospectral set of isometry classes of flat metrics
  on $\Sigma$ is compact in the $C^\infty$-topology.
\end{theorem}
Let $\{\rho_i \}_{i=1}^\infty$ be a sequence of flat (i.e.
sectional curvature is zero) metrics on $\Sigma$ having the same
Dirichlet spectrum. By heat kernel asymptotic expansions (see \cite{MS}),
all such metrics $\rho_i$'s have the same area, say $A$, as
well as the same boundary length, say $L$.  Write $\rho_i =
e^{2\phi_i} \sigma_i$ where  $\sigma_i$ is a uniform metric of
type II of area $A$. Each $\sigma_i$ induces an element $[\sigma_i ] \in
\mathcal{M}_{II} (\Sigma, A)$.
\begin{proposition}\label{compactness of sigma}
  The above sequence $\{ [\sigma_i] \}$ is compact in
  $\mathcal{M}_{II} (\Sigma, A)$.
\end{proposition}
\begin{proof}
  We drop the index $i$ for a moment. Consider
  $\tau = \Psi (\sigma)$, i.e.  $\tau$ is a uniform metric of type I
  of area $A$, and write $\rho= e^{2\psi} \tau$.
  Without loss of generality we may assume $A= - 2\pi \chi(\Sigma)$
  and so the curvature $K_\tau$ of $\tau$ is $-1$.
  Polyakov-Alvarez formula reads

  \begin{align}\label{for compactness of sigma}
    h(e^{2\psi}\tau) = &\frac{1}{6 \pi}
    \big{\{} \frac{1}{2}\int_\Sigma | \nabla_\tau \psi |^2 \,dA_\tau
    -\int_\Sigma  \psi \, dA_\tau  \big{\}}
    + \frac{1}{4 \pi} \int_{\partial \Sigma} \partial_n ^\tau \psi \,ds_\tau
    + h(\tau)
  \end{align}
  Since by \eqref{conformal curvature}
  $ \Delta_\tau \psi   =1 $, we see
  \begin{align*}
    \int_{\partial \Sigma} \partial_n ^\tau \psi ds_\tau
    = -\int_\Sigma \Delta_\tau \psi dA_\tau
    = 2 \pi \chi(\Sigma)
  \end{align*}
  and by Jensen's inequality
  \begin{align*}
    \int_\Sigma  2 \psi  dA_\tau  /A \leq \log \big{(}\int_\Sigma
    e^{2\psi} dA_\tau / A \big{)}  = \log 1 = 0.
  \end{align*}
  So back to \eqref{for compactness of sigma} we get
  \begin{align}\label{h(tau) ineq}
    h(e^{2\psi}\tau) \geq  \frac{\chi(\Sigma)}{2}    + h(\tau).
  \end{align}
  For each $\rho_i$,
  consider the corresponding unform metric $\tau_i$
  of type $I$, i.e. $\tau_i = \Psi (\sigma_i)$.
  By isospectrality,
  $h(\rho_i) = \text{const.} $ for all $i$. Thus by \eqref{h(tau) ineq},
  $h(\tau_i) \leq \text{const}$.
  From Theorem~\ref{T:proper on moduli I}, $h$ is proper and bounded below
  on $\mathcal{M}_{I}(\Sigma, A)$, and so
  the sequence $\{ [\tau_i ]\}$ induced in $\mathcal{M}_{I}(\Sigma, A)$
  is compact. Apply the map $\tilde{\Phi}$
  (see Theorem~\ref{T:tilde Psi}), to see $\{ [\sigma_i ] \}$ is compact.
\end{proof}
\begin{remark}
    It seems hard to get directly the compactness of
    the type II uniformization of isospectral flat metrics,
    without working on the type I uniformization first then using
    Theorem~\ref{Psi} as we did in the above proof.
    The technical difficulty exists because the corresponding
    Polyakov-Alvarez  formula has the term involving boundary geodesic
    curvature of the
    type II uniform metrics, which we do not have good control of when the area is fixed.
    This is true especially
    if the metrics are near the boundary of the corresponding moduli space,
    a case \emph{a priorily} possible though excluded by
    Proposition~\ref{compactness of sigma}.
\end{remark}
From Proposition~\ref{compactness of sigma}, we can find isometric representatives of
$\sigma_i $'s which consist a compact set in $C^\infty$,
and we denote these representatives by the same $\sigma_i$'s.
The remaining part is in principle the same as
in the last  section of \cite{OPS3}. It is included  for reader's convenience.
Note that  by compactness of $\{ \sigma_i \}$,
all the Sobolev $t$-norms $\| \cdot \|_t$ of functions on $\Sigma$,
resp. $|\cdot |_t$ of functions on $\partial \Sigma$,
resp. the $C^j$ norms, induced by $\sigma_i$ are
equivalent by uniformly bounded constant multiples;
therefore, in the consideration below, we may
without loss of generality deal with only such norms for one fixed metric.
Note also that  all the quantities induced by the metric $\sigma_i$
are uniformly bounded
in the $C^j$ norm for each $j$.
Now, for compactness of the (isometric representative) sequence $\{ \rho_i \}$,
it is enough to show that $\{\phi_i \}$ has
a convergent subsequence in $C^\infty$. By
a diagonal argument using Rellich's compactness it is equivalent
to show that the sequence $\{\phi^i\}$ is uniformly bounded in
$\| \cdot\|_t $ for each integer $t$.

From this point on, we will drop the index $i$.
Using \eqref{conformal curvature} and\eqref{geodesic curvature},  we see
\begin{align}\label{harmonicity}
  \Delta_\sigma \phi = 0  , \\
  \label{geodesic curvature for compactness}
  k =e^{-\phi} ( k_\sigma  + \partial_n ^\sigma \phi )
\end{align}
where $k_\sigma =\frac{2\pi \chi(\Sigma)}{\int_{\partial \Sigma} ds_\sigma}$
and $\partial_n ^\sigma$ is now  the Dirichlet-Neumann operator.
By Lemma~\ref{elliptic harmonic}, it is enough to show
that the $\phi$'s are uniformly bounded for $| \cdot |_t$
for each integer $t$. We will do an induction on $t$.

\noindent {\bf Step 1}\cite{OPS3}
First let's deal with the Sobolev $1$-norm of $\phi$.
The Polyakov-Alvarez formula reads
\begin{align}\label{Polyakov for bounds}
  h(e^{2\phi} \sigma) = \frac{1}{6 \pi}
  \big{\{} \frac{1}{2}\int_\Sigma | \nabla_\sigma \phi |^2 \,dA_\sigma +
  k_\sigma \int_{\partial \Sigma}  \phi ds_\sigma \big{\}}+  h(\sigma).
\end{align}
By isospectrality,
\begin{align}\label{by isospectrality}
  h(e^{2\phi} \sigma) = \text{const.},
\end{align}
and $h(\sigma)$ is bounded above and below by the compactness of
the set $\{\sigma_i \}$.
Since $\int_{\partial \Sigma} e^\phi ds_\sigma = L$,
by Jensen's inequality
\begin{align*}
  \int_{\partial \Sigma} \phi ds_\sigma \leq
  \log (\frac{L}{\int_{\partial \Sigma} ds_\sigma })
  \int_{\partial \Sigma} ds_\sigma  \leq \text{const.}.
\end{align*}
Combining this with \eqref{Polyakov for bounds},
\eqref{by isospectrality}, and the fact that $k_\sigma < 0$ and
$|k_\sigma |$ is bounded, we see
\begin{align}\label{by Jensen 2}
  | k_\sigma \int_{\partial \Sigma} \phi ds_\sigma | \leq  \text{const}.
\end{align}
and
\begin{align}\label{bound int phi}
  |  \int_{\partial \Sigma} \phi ds_\sigma |   \leq \text{const.}.
\end{align}

By \eqref{Polyakov for bounds}, \eqref{by isospectrality},
and \eqref{by Jensen 2},
\begin{align}\label{grad bound}
  \int_\Sigma | \nabla_\sigma \phi |^2 dA_\sigma \leq \text{const.}.
\end{align}

By \eqref{bound int phi},
\eqref{grad bound}, \eqref{1/2 sobolev}, and \eqref{Trudinger inequality},
\begin{align}
  \label{1/2 bound phi}
  |\phi|_{1/2}  \leq \text{const.},
  \\
  \label{bound of exp 2phi}
  \int_{\partial \Sigma} e^\phi ds_\sigma , \ \
  \int_{\partial \Sigma} e^{2 \phi} ds_\sigma \leq \text{const.} .
\end{align}

By Melrose's result \cite{Me} for an isospectral set of flat metrics,
 the geodesic curvature $k$
as a function on $\partial \Sigma$ is
uniformly bounded in the $C^j$ norm for each $j$.
(His original result was in Euclidean context
but it can be easily carried over to a flat surface case.)
Since $k$ and $k_\sigma$ are uniformly bounded, we see
by \eqref{geodesic curvature for compactness},
\begin{align*}
  |  \partial_n ^\sigma \phi    | ^2
  \leq \text{const.} ( e^{2 \phi}  + e^\phi + 1),
\end{align*}
and  \eqref{bound of exp 2phi} gives
\begin{align}\label{normal derivative bound}
  \int_{\partial \Sigma} |  \partial_n ^\sigma \phi    | ^2   ds_\sigma
  & \leq \text{const.} ( \int_{\partial \Sigma}
  e^{2 \phi} ds_\sigma +  \int_{\partial \Sigma} e^\phi ds_\sigma  +
  \int_{\partial \Sigma}  ds_\sigma )
  \leq \text{const.}.
\end{align}

\begin{lemma}\label{Sobolev norm}\emph{(See \cite{OPS3}: eqn. (5.3))}
    For a smooth function $f$ on $\partial \Sigma$,
  \begin{align*}
  |\partial_n f |_0 + | f|_0 = |f |_1  ,
  \end{align*}
  where $\partial_n$ is the Dirichlet-Neumann operator.
\end{lemma}
So by \eqref{1/2 bound phi} and \eqref{normal derivative bound}
we get the Sobolev $1$-norm $|\cdot |_t$ bound for $\phi$.
In particular, the $\phi$'s are uniformly bounded on $\partial \Sigma$.

\noindent {\bf Step 2}\cite{OPS3}
From \eqref{geodesic curvature for compactness},
for the derivative $\partial_s$ along $\partial \Sigma$
\begin{align*}
  | \partial_s ^{ \, t } \partial_n \phi |_0
  = |\partial _s ^{ \, t} (-k_\sigma + e^\phi k )|_0
  \leq \hbox{const.}_t \  | \phi |_t        .
\end{align*}
The Dirichlet-Neumann operator $\partial_n$ is
a pseudodifferential operator of order $1$ and
\begin{align*}
  | [ \partial_s ^{\, t} , \partial_n ] \phi |_0 \leq
  \hbox{const.}_t \ |\phi |_t ,
\end{align*}
where $[\ ,\ ]$ is the commutator operation.
Now by Lemma~\ref{Sobolev norm},
\begin{align*}
  |\phi|_{t+1} \leq | \partial_n \partial_s ^{\, t} \phi |_0
  + | \phi |_t \leq \hbox{const.}_t \ |\phi|_t   .
\end{align*}
This completes the induction and finishes the proof
of the $C^\infty$-compactness of the set $\{ \phi_i\}$
and so of $\{\rho_i \}$, and so finalizes the proof of
Theorem~\ref{T:isospectral compactness}.

\section{Properness of the height on $\mathcal{M}_{II}(\Sigma, A)$}\label{S:flat proper}
Let $\sigma \in \mathbf{M}_{II} (\Sigma, A)$ and
$\tau \in \mathbf{M}_I (\Sigma, A)$ such that
$\sigma = e^{2\psi} \tau$.
Proceed exactly as for \eqref{h(tau) ineq} to get
\begin{align}\label{h(sigma) h(tau) ineq}
  h(e^{2\psi}\tau) \geq   \frac{1}{2}\chi(\Sigma)     + h(\tau).
\end{align}
This allows us to show the following theorem.
\begin{theorem}\label{T:proper on moduli II}
  For each $A >0$,
  the height $h$ is a proper function on
  the moduli space $ \mathcal{M}_{II} (\Sigma, A)$, i.e.
  \begin{align*}
    h(M) \rightarrow + \infty
  \end{align*}
  as the isometry class $[M]$ approaches
  $ \partial \mathcal{M}_{II}(\Sigma, A) $.
\end{theorem}
\begin{proof}
  The theorem follows from \eqref{h(sigma) h(tau) ineq},
  Theorem~\ref{T:tilde Psi}, and Theorem~\ref{T:proper on moduli I}.
\end{proof}

\begin{remark}
    It is easy to see that for two conformally equivalent uniform
    metrics $\tau$ and $\sigma$, respectively of type I and II, of
    the same boundary length,
    \begin{align*}
        h(\tau) \ge \frac{1}{2}\chi(\Sigma) + h (\sigma).
    \end{align*}
\end{remark}

\section{Further Remarks}\label{S:remarks}
It is desirable to find a connection between our results and the results of
Khuri \cite{Kh}, Osgood, Phillips, and Sarnak \cite{OPS3}: they showed the properness (resp. non-properness) of height function \cite{OPS3} (resp. \cite{Kh})
for $(0, n)$-type surfaces, $n >0$ (resp.  for $(g, n)$-type surfaces, $g>0, n>0$),
when the boundary length is fixed instead of the area.
Note that
the two conditions, one fixes the area
and the other fixes the boundary length, are equivalent by scaling,
and there is
a nice formula \eqref{scale height} for the change of height
under scaling.

It seems interesting to get a better understanding of the maps
$\tilde{\Psi}$ and $\tilde{\Phi}$ between the two spaces
$\mathcal{M}_{I}(\Sigma, A)$ and $\mathcal{M}_{II}(\Sigma, A)$. For
example, it would be nice to investigate the extension of these maps
to the compactifications of $\mathcal{M}_I(\Sigma, A)$ and
$\mathcal{M}_{II}(\Sigma, A)$, and it would also be interesting to
see the geometric properties of these maps with respect to certain
natural metrics given on those moduli spaces.

Finally, we wonder whether the set of isospectral Riemannian metrics
(without any further restriction) on a compact surface with boundary
(or a higher dimensional manifold) is compact in the $C^\infty$-topology. One may try to modify the method
of Osgood, Phillips, and Sarnak in \cite{OPS2}.

\bibliographystyle{plain}

\end{document}